\RequirePackage{fix-cm}
\documentclass[smallextended]{svjour3}       
\smartqed  
\usepackage{url}

\usepackage{amsmath,amsfonts,amssymb,graphicx}
\usepackage{tikz}
\usepackage{multirow}
\usepackage[normalem]{ulem}
\usepackage{comment}

\renewcommand{\b}[1]{\boldsymbol{#1}}
\newcommand{\bc}{\b{c}}
\newcommand{\bcp}{\b{c'}}
\newcommand{\bctp}{\b{\tilde c'}}

\newcommand{\cT}{\mathcal{T}}

\newcommand{\ddiv}{\operatorname{div}}

\newcommand{\hP}{\widehat{P}} 
\newcommand{\hPi}{\widehat{\Pi}}
\newcommand{\hE}{\widehat{E}} 
\newcommand{\Hdiv}[1][\Omega]{H(\ddiv,#1)}

\newcommand{\linspan}{\operatorname{span}}
\newcommand{\norm}[1]{\|#1\|}

\newcommand{\R}{\mathbb{R}}

\newcommand{\IndEff}[1]{{\mathcal{I}( #1 )}} 

\journalname{Numerische Mathematik}
\begin{document}

\title{Fully computable a posteriori error bounds for eigenfunctions
\thanks{
  The first author is supported by Japan Society for the Promotion of Science: Fund for the Promotion of Joint International Research (Fostering Joint International Research (A)) 20KK0306, 
  Grant-in-Aid for Scientific Research (B) 20H01820, 21H00998, and Grant-in-Aid for Scientific Research (C) 18K03411.
  The second author is supported by the Czech Science Foundation, project no.~20-01074S, and by RVO 67985840. 
  This work is also supported by the Research Institute for Mathematical Sciences, an International Joint Usage/Research Center of Kyoto University.}
}


\author{Xuefeng Liu \and Tom\'a\v s Vejchodsk\'y }


\institute{Xuefeng Liu \at
              Graduate School of Science and Technology, Niigata University, 8050 Ikarashi 2-no-cho, Nishi-ku, Niigata City, Niigata 950-2181, Japan\\
              \email{xfliu@math.sc.niigata-u.ac.jp}           
           \and
           Tom\'a\v s Vejchodsk\'y \at
           Institute of Mathematics, Czech Academy of Sciences, \v{Z}itn\'a 25, Prague 1, 115\,67, Czech Republic \\
           \email{vejchod@math.cas.cz} 
}

\date{Received: date / Accepted: date}

\maketitle

\begin{abstract}
  For compact self-adjoint operators in Hilbert spaces, two algorithms are proposed to provide fully computable \emph{a posteriori} error estimate for eigenfunction approximation. 
   Both algorithms apply well to the case of tight clusters and multiple eigenvalues, under the settings of target eigenvalue problems. 
  Algorithm I is based on the Rayleigh quotient and the min-max principle that characterizes the eigenvalue problems. The formula for the error estimate provided by Algorithm I is easy to compute and applies to problems with limited information of Rayleigh quotients. 
  Algorithm II,
as an extension of  the Davis--Kahan method, 
takes advantage of the
  dual formulation of differential operators along with the Prager--Synge technique and provides greatly improved accuracy of the estimate, especially for the finite element 
  approximations of eigenfunctions.
  Numerical examples of eigenvalue problems of matrices and the Laplace operators over convex and non-convex domains  illustrate the efficiency of the proposed algorithms.
  

\end{abstract}

\section{Introduction}

The weak form of eigenvalue problems for linear elliptic partial differential operators motivates us to consider an abstract setting of a pair of Hilbert spaces
$V$ and $W$ with inner products $a(\cdot, \cdot)$ and $b(\cdot, \cdot)$, respectively, and a compact operator $\gamma: V\to W$ between them. We seek eigenvalues $0 < \lambda_1 \leq \lambda_2 \leq \cdots$, repeated according to their multiplicity, and corresponding eigenfunctions $u_i \in V\setminus\{0\}$,  $i=1,2,\dots$ such that
\begin{equation}
  \label{eq:eigp}
  a(u_i, v) = \lambda_i b(\gamma u_i, \gamma v)\quad \forall v \in V.
\end{equation}
Section~\ref{se:eigenproblem} provides more details about the well-posedness of this problem and gives examples how the (generalized) matrix, Laplace, and Steklov eigenvalue problems fit into this setting.

The problem to determine eigenvalues $\lambda_i$ is well posed in the sense that small perturbations of the data lead to small perturbations of eigenvalues. 
However, the variation of eigenfunctions $u_i$ 
upon a small perturbation of the data is not necessarily small, and can even be discontinuous. For example, if two close and simple eigenvalues merge to one multiple eigenvalue then the two corresponding orthogonal eigenfunctions abruptly change into a two dimensional eigenspace. Thus, eigenfunction determination in case of tightly clustered or multiple eigenvalues is an ill-conditioned problem. 
Any attempt to estimate the error of approximate eigenfunctions has to take into the account this ill-conditioning.

Our approach is to consider the space spanned by eigenfunctions corresponding to all eigenvalues within a cluster.
This space is well conditioned provided the cluster is well separated from the rest of the spectrum.
We propose error estimators that bound the \emph{directed distance} \cite[\S5.15]{Meyer:2000} between the approximate and the exact space spaces of eigenfunctions in norms induced by both inner products $a$ and $b$, see estimates \eqref{eq:Deltaest} and \eqref{eq:Deltaest_b} below.
Further, we present a bound on the distance in the $a$-sense obtained from the distance in the $b$-sense; see \eqref{eq:a_by_b} below.

The error estimators in our approach provide rigorous upper bounds on these distances without any {\em a priori} information about the approximate eigenfunctions.
These estimators are fully computable in terms of approximate eigenfunctions and two-sided bounds on eigenvalues 
and their quality depends on the width of clusters and spectral gaps between them.
For the Laplace eigenvalue problem, the proposed estimators generalize the idea from \cite{BirBooSwaWen:1966}; see Remark \ref{rem:birkhoff-etl-result} below.

To compute the needed two-sided bounds on individual eigenvalues,
we use the recently developed approach based on the finite element method with 
explicit error estimation
\cite{Liu2015,Liu2020} (see also, \cite{LiuOis2013,CarGal2014,CarGed2014}) for the lower bounds on eigenvalues
and the Lehmann--Goerisch method \cite{Lehmann1949,Lehmann1950,GoeHau1985} 
for their high-precision improvements. 
Note that the Lehmann--Goerisch method should be attributed to T. Kato as well, because his independently developed method \cite{Kato1949}, gives essentially the same bounds as Lehmann's method.
In the current paper, we focus on the estimation of eigenfunctions and the two-sided bounds of eigenvalues are assumed to be known.

Error estimates for symmetric elliptic partial differential eigenvalue problems are widely studied in the literature. We refer to classical works \cite{Chatelin1983,BabOsb:1991,Boffi:2010} for the fundamental theories. 
The majority of the existing literature concerns error estimates valid asymptotically or containing unknown constants; see, e.g., \cite{ArmDur2004,DarDurPad2012,DurGasPad1999,giani2018posteriori,GiaHal2012,HuHuaLin2014,JiaCheXie2013,MehMie2011,Yang2010}.
Recently, fully computable (containing no unknown constants) and guaranteed (bounding the error from above on all meshes, not only asymptotically) error estimates for eigenvalue problems appeared. 
Papers  \cite{CarGal2014,CarGed2014,Liu2015,Liu2020,LiuOis2013,SebVej2014,Vejchodsky2018b,Vejchodsky2018} concern the eigenvalues.
Particularly, as a general framework, the method proposed in \cite{Liu2015,Liu2020} has been applied to eigenvalue problems of various differential operators, including the Stokes operator \cite{Xie2LIU-2018}, 
the Steklov operator \cite{you-xie-liu-2019}, and biharmonic operators related to the quadratic interpolation error constants \cite{liu-you:2018,LiaoYuLiu2019}.

\medskip

Concerning eigenfunctions, 
Davis and Kahan \cite{davis1970rotation} provide a fundamental result
for the estimation of the distance of eigenspaces utilizing the strong residual $R=A\hat{x}-\hat{\lambda} \hat{x}$, where $\{ \hat{x}, \hat{\lambda}\}$ stands for an approximate eigenpair of a matrix (or operator) $A$. 
Recent paper \cite{Nakatsukasa2020} extends Davis--Kahan's approach by further orthogonal decomposition of $A$ with respect to the trial eigenspaces.

Similar to Davis--Kahan's approach, 
papers \cite{CanDusMadStaVoh2017,CanDusMadStaVoh2018,CanDusMadStaVoh2019,CarGed2014,HongXieYueZhang2018} 
consider the weak residual. These methods provide guaranteed, robust, and optimally convergent {\em a posteriori} bounds for eigenvalues and corresponding eigenfunctions for finite element approximations assuming an \emph{a priori} knowledge of bounds on eigenvalues. All these existing approaches introduce dual functions from $H(\mbox{div};\Omega)$ to approximate the gradient of eigenfuctions and need to solve an auxiliary problem to find them. 
Alternative approach \cite{toyonaga2002verified} formulates the error estimation problem as a fixed-point problem.
{Section 6 contains the comparison of the 
efficiency of our proposed method and the one in \cite{CanDusMadStaVoh2019}.}

\medskip

In contrast to the existing literature, we propose two algorithms to  estimate error of approximate eigenfunctions. Proposed error estimates have the following features.

\begin{itemize}
\item  Algorithm I in \S \ref{se:ebound} is based on the Rayleigh quotient and the min-max principle that characterizes the eigenvalue problem. This algorithm only utilizes Rayleigh quotients of approximate eigenfunctions and thus the estimator is easy to compute, especially when compared to Algorithm II and the existing literature. 
A defect of this algorithm is that the estimate for the concerned cluster depends on previous clusters {and the} width of the concerned cluster.

\item Algorithm II in \S \ref{se:est-residual-error} utilizes the residual error for variational representation of differential operators.
Here, Davis--Kahan's method, originally proposed for strongly formulated eigenvalue problems, is successfully extended to weakly formulated eigenvalue problems. 
In case of the Laplace eigenvalue problem, we obtain sharper estimate for the error of approximate eigenfunctions by further utilizing an auxiliary flux $p_h \in H(\mbox{div})$.

The estimates by Algorithm II are independent from the cluster index, its width (difference of the largest and smallest eigenvalue in the cluster), and its degree (number of eigenvalues in the cluster). 
The quality of the estimate depends solely on the residual error and the spectral gap (distance of the cluster from its neighboring clusters). Numerical examples presented below show that this estimate is very accurate, 
especially for finite element approximations.
\end{itemize}

\medskip

The rest of the paper is organized as follows.
Section~\ref{se:eigenproblem} briefly recalls the main properties of the abstract eigenvalue problem.
Section~\ref{se:ebound} presents the {\em a posteriori} error bounds for eigenfunctions in both $a$- and $b$-sense.
Section~\ref{se:a_by_b} derives an error bound for the directed distance measured in the $a$-sense based on the error bound in the $b$-sense. 
Section~\ref{se:est-residual-error} utilizes the residual error of the differential operators and extended Davis--Kahan's theorem to obtain the error estimation.
Section~\ref{se:numex} provides numerical examples for matrix and Laplace eigenvalue problems.
Finally, Section~\ref{se:conclusions} draws the conclusions.

\section{Eigenvalue problem for a compact self-adjoint operator}
\label{se:eigenproblem}

As we mentioned above, weak formulations of elliptic partial differential eigenvalue problems motivate us to consider Hilbert spaces $V$ and $W$ with inner products $a(\cdot, \cdot)$ and $b(\cdot, \cdot)$, respectively, and a compact operator $\gamma: V\to W$ between them.
Let $T:W\to V$ be the solution operator such that for given $f\in W$, the element $Tf \in V$ is uniquely determined by the identity
\begin{equation}
\label{eq:def-T-operator}
a(Tf, v) = b(f,\gamma v)\quad \forall v \in V.
\end{equation}
It is easy to see that $T \gamma : V\to V$ is a compact self-adjoint operator. The eigenvalue problem for the operator $T \gamma$ is equivalent to finding eigenfunctions $u_i \in V\setminus\{0\}$ and corresponding eigenvalues $\mu_i \ge 0 $ such that 
\begin{equation}
  \label{eq-main-eig-form}
  b(\gamma u_i, \gamma v) = \mu_i \, a(u_i, v)  \quad \forall v \in V.
\end{equation}

The Hilbert--Schmidt Theorem (see, e.g., \cite[Theorem 4.10.1]{debnath2005introduction}) 
provides the existence of a countable sequence of (finite) non-zero eigenvalues
$$
 0\le  \lim_{i\to d} \mu_i \leq \cdots \leq \mu_2 \leq \mu_1, 
$$ 
of the operator  $T \gamma$, where $d = \operatorname{dim} V$ can be infinite. 
Note that these eigenvalues are repeated according to their multiplicities and  
that if $\mbox{dim}(\operatorname{Ker}(T \gamma))>0$ then there are also zero eigenvalues of $T \gamma$.

Eigenvalue problems for linear elliptic partial differential operators are usually formulated for the reciprocals of $\mu_i$. Denoting $\lambda_i = 1/\mu_i$ for $\mu_i>0$, we clearly have
$$
  0 < \lambda_1 \leq \lambda_2 \leq \cdots,
  \quad\text{and}\quad 
  a(u_i, v) = \lambda_i b(\gamma u_i, \gamma v)\quad \forall v \in V,
$$
namely the eigenvalue problem \eqref{eq:eigp}. 

The eigenfunctions $u_i \in V$ are assumed to be normalized such that
$$
a(u_i,u_j) = \delta_{ij}, \quad i,j = 1,2, \dots,
$$
where $\delta_{ij}$ stands for the Kronecker delta.
With this normalization, the Hilbert--Schmidt Theorem also implies that
any $v \in V$ can be expressed as
\begin{equation}
  \label{eq:Parseval2}
  v =  \sum_{i=1}^d a(v,u_i)  u_i +v_0, \quad\text{where } v_0 \in \operatorname{Ker}(T \gamma).
\end{equation}
Noticing that $b(\gamma v_0, \gamma  v)= a(T \gamma v_0, v)=0$ for all $v\in V$, we have
\begin{equation}
\label{eq:Parseval}
\norm{v}_a^2 =  \sum_{i=1}^d |a(v,u_i)|^2 + \|v_0\|_a^2 
\quad\text{and}\quad
\norm{\gamma v}_b^2 =  \sum_{i=1}^d |a(v,u_i)|^2 /\lambda_i
\quad\forall v \in V.
\end{equation}

\begin{remark}\label{re:examples}
The Laplace eigenvalue problem in a domain $\Omega$ with homogeneous Dirichlet boundary conditions fits the above setting for
\begin{equation*}
  \left\{
    \begin{array}{l} 
      \displaystyle V = H_0^1(\Omega), \quad a(u,v)=\int_{\Omega} \nabla u \cdot \nabla v \,\mbox{d}x; \\
      \displaystyle W=L^2(\Omega), \quad b(u,v)=\int_\Omega uv \,\mbox{d}x; \\
    \gamma \mbox{ being the identity operator}.  
    \end{array}
  \right.
\end{equation*}
The Steklov eigenvalue problem fits this settings for
\begin{equation*}
  \left\{
    \begin{array}{l} 
      \displaystyle V=H^1(\Omega), \quad a(u,v)=\int_{\Omega} \nabla u \cdot \nabla v + uv \,\mbox{d}x;\\
      \displaystyle W=L^2(\partial \Omega), \quad b(u,v)=\int_{\partial \Omega} uv \,\mbox{d}s;\\
      \gamma \mbox{ being the trace operator}.  
    \end{array}
  \right.
\end{equation*}
The generalized matrix eigenvalue problem $A u_i = \lambda_i B u_i$ for a symmetric positive definite matrix $A \in \R^{d\times d}$ and a symmetric positive semidefinite matrix $B \in \R^{d\times d}$ fits the above setting for $d<\infty$, $V = \R^d$, $W = \operatorname{range} B$, $\gamma = B^{1/2}$, $a(u,v) = v^T A u$, and $b(u,v) = v^T u$. 
\end{remark}

\begin{figure}[ht]
  \begin{tikzpicture}[scale=1]
  \newcommand{\tlen}{0.1}
  \newcommand{\tick}[1]{\draw [semithick] (#1,-\tlen)--(#1,\tlen);}
  \draw [semithick] (0,0)--(7.5,0);
  \draw [semithick,dotted] (7.5+0.2,0)--(7.5+1.3,0);
  \draw [semithick] (7.5+1.5,0)--(7.5+5,0);
  
  \tick{0.5}\node [below] at (0.5,-\tlen) {$0$};
  
  \tick{1.7}\node [above] at (1.7,\tlen) {$\lambda_{n_1}$};
  \tick{1.9}
  \tick{2}
  \tick{2.2}
  \tick{2.4}
  \tick{2.6}
  \tick{2.7}
  \tick{3}\node [above] at (3,\tlen) {$\lambda_{N_1}$};
  
  \tick{5}\node [above] at (5,\tlen) {$\lambda_{n_2}$};
  \tick{5.25}
  \tick{5.4}
  \tick{5.5}
  \tick{5.8}
  \tick{6}
  \tick{6.2}
  \tick{6.35}
  \tick{6.5}
  \tick{6.7}
  \tick{6.8}
  \tick{7}\node [above] at (7,\tlen) {$\lambda_{N_2}$};
  
  \tick{9.5}\node [above] at (9.5,\tlen) {$\lambda_{n_K}$};
  \node [below] at (9.5,-\tlen) {$\lambda_n$};
  \tick{9.65}
  \tick{9.75}
  \tick{10}
  \tick{10.1}
  \tick{10.35}
  \tick{10.5}
  \tick{10.7}
  \tick{10.8}
  \tick{11}\node [above] at (11,\tlen) {$\lambda_{N_K}$};
  \node [below] at (11,-\tlen) {$\lambda_N$};
  \end{tikzpicture}
  \caption{Clusters of eigenvalues on the real axis.}
  \label{fi:clusters}
  \end{figure}
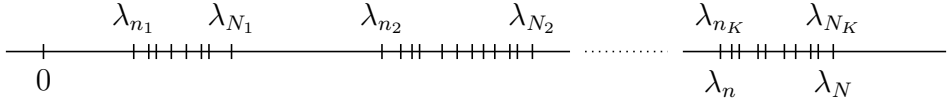

In order to formulate the bound on eigenfunctions, a notation for clusters of eigenvalues has to be introduced.
Let us focus on the leading $K$ clusters.
Let $n_k$ and $N_k$ stand for indices of the first and the last eigenvalue in the $k$-th cluster, $k=1,2,\dots,K$, respectively. In particular, $n_1 = 1$, $n_{k+1} = N_k + 1$, and the $k$-th cluster is formed of $N_k - n_k + 1$ eigenvalues $\lambda_{n_k}$, $\lambda_{n_k + 1}$, \dots, $\lambda_{N_k}$; see Figure~\ref{fi:clusters}. Notice that the eigenvalues in a cluster do not necessarily equal to each other.
To simplify the notation, we set $n=n_K$ and $N=N_K$.

Each cluster is associated with the space
$E_k = \linspan\{ u_{n_k}, u_{n_k+1}, \dots, u_{N_k} \}$ of exact eigenfunctions.
Similarly, arbitrary \emph{linearly independent} approximations $\hat u_i \in V$ of exact eigenfunctions $u_i$, $i=1,2,\dots,N_K$, form the corresponding approximate spaces
$\hE_k = \linspan\{ \hat u_{n_k}, \hat u_{n_k+1}, \dots, \hat u_{N_k} \}$. Thus $\mbox{dim}(\hE_k)=N_k-n_k+1$.
Spaces $\hE_k$, $k=1,2,\dots,K$, of approximate eigenfunctions need not be orthogonal to each other. 

\begin{remark}
Given an eigenvalue problem, the {\em a priori} information about the distribution of 
eigenvalues is usually unknown. Therefore, in practical problems, the partition of the spectrum into clusters is determined numerically using two-sided bounds on eigenvalues. 

If the enclosing intervals of eigenvalues overlap, it is natural to define a cluster for such eigenvalues. 
Consequently, the internal distribution of eigenvalues inside the cluster is regarded as  uncertain and only the lower and upper bound on the whole cluster will be used for analysis. 

\end{remark}

\section{Algorithm I: Estimation for eigenfunctions based on the Rayleigh quotients}
\label{se:ebound}

The goal of this section is to derive an estimate of the directed distance between spaces $E_K$ and $\hE_K$ of exact and approximate eigenfunctions for the
$K$-th cluster. The directed distance between two subspaces $E$ and $\hE$ of $V$ is defined both in the $a$- and $b$-sense as follows; see e.g.,  \cite[\S5.15]{Meyer:2000}, 
\begin{equation}
\label{eq:Delta}
\delta_a(E, \hE) = \max_{\substack{v \in E\\ \norm{v}_a=1}} \min_{ \hat{v} \in \hE} \| v - \hat{v} \|_a
\quad\text{and}\quad
\delta_b(E, \hE) = \max_{\substack{v \in E\\ \norm{\gamma v}_b=1}} \min_{ \hat{v} \in \hE} \| \gamma (v - \hat{v}) \|_b.
\end{equation}
The directed distance is not symmetric, in general, but
if the two subspaces are of the same finite dimension, which is the case in this paper, 
then the directed distance is symmetric and coincides with the gap between spaces.
The following characterization of the directed distance will be used in the proof of the main theorem; see also \cite[\S5.15]{Meyer:2000}. 

Note that in \cite{CanDusMadStaVoh2019},  the Hilbert–Schmidt norm for the cluster residual error is utilized to estimate the error of approximate eigenfunction. 

\begin{lemma}\label{le:delta}
  Let $E$ and $\hE$ be the two finite dimensional subspaces of the Hilbert space $V$.
  Let $\hP:V\to \hE$ be the orthogonal projector with respect to $a(\cdot, \cdot)$. 
Then, the directed distance defined in \eqref{eq:Delta} satisfies
\begin{alignat}{2}
  \label{eq:Delta2}
  \delta_a(E,\hE) &= 
  \max_{\substack{v \in E \\ \norm{v}_a = 1}} \norm{v - \hP v}_a,
&\qquad
  \delta_b(E,\hE) &= 
  \max_{\substack{v \in E \\ \norm{\gamma v}_b = 1}} \norm{\gamma(v - \hP v)}_b,
\\
  \label{eq:Delta3}
  \delta_a^2(E,\hE)
  &= 1 - \min_{\substack{v \in E\\ \norm{v}_a = 1}} \norm{\hP v }_a^2,
  &\qquad
  \delta_b^2(E,\hE)
  &= 1 - \min_{\substack{v \in E\\ \norm{\gamma v}_b = 1}} \norm{\gamma \hP v }_b^2.
\end{alignat}
\end{lemma}
\begin{proof} 
  The first statement in \eqref{eq:Delta2} follows immediately from \eqref{eq:Delta}, because the orthogonal projection $\hP v$ is the closest element to $v$ in the subspace $\hE$. 
  For the second statement, we apply the same idea to the orthogonal projector $\hPi: \gamma V \to \gamma\hE$ with respect to the scalar product $b(\cdot,\cdot)$. This projector satisfies $\hPi \gamma = \gamma \hP$ and, hence, the second statement holds true.
  Identities \eqref{eq:Delta3} follow from \eqref{eq:Delta2} due to the equality $\norm{\hP v}_a^2 + \norm{v - \hP v}_a^2 = \norm{v}_a^2$ and a similar equality for $\hPi$.
\end{proof}

\begin{remark}
For one dimensional spaces $E=\linspan\{u\}$ and $\hE = \linspan\{\hat{u}\}$, we have $\delta_a(E,\hE) = \sin \alpha$, where $\alpha$ is the angle between $u$ and $\hat{u}$. Consequently, the difference of $u$ and $\hat{u}$ in the $a$-norm can be expressed by the cosine theorem as
$$
  \norm{u - \hat{u}}_a^2 = \norm{u}_a^2 + \norm{\hat{u}}_a^2 - 2 \norm{u}_a \norm{\hat{u}}_a \sqrt{1 - \delta_a^2(E,\hE)}.
$$
Hence, an upper bound on $\delta_a(E,\hE)$ yields immediately an upper bound on the error $\norm{u - \hat{u}}_a$. Moreover, if $\norm{u}_a =\norm{\hat{u}}_a = 1$ then $\norm{u - \hat{u}}_a^2 = \delta_a^2(E,\hE) + O(\delta_a^4(E,\hE))$ and the difference between the directed distance $\delta_a(E,\hE)$ and the $a$-norm of the error is negligible for small $\delta_a(E,\hE)$.
The same conclusions clearly hold also for the norm induced by $b$.
\end{remark}

The main result of this section is presented in Theorem~\ref{th:mainenergy} below.
In order to formulate it, we introduce 
\emph{measures of non-orthogonality}
between finite dimensional subspaces $E$ and $E'$ of $V$ both in the $a$- and $b$-sense as
\begin{equation}
  \label{eq:defzeta}
  \hat{\varepsilon}_a(E,E') = \max_{\substack{v\in E\\ \| v\|_a=1}} \max_{\substack{v'\in E'\\ \| v'\|_a=1}} a(v, v')
\quad\text{and}\quad  
  \hat{\varepsilon}_b(E,E') = \max_{\substack{v\in E\\ \|\gamma v\|_b=1}} \max_{\substack{v'\in E'\\ \|\gamma v'\|_b=1}} b(\gamma  v , \gamma v' ).
\end{equation}
Both these measures of non-orthogonality can be easily computed or estimated by using the following lemma. To simplify the exposition, we present the result for $\hat{\varepsilon}_a$ only, because the statement for $\hat{\varepsilon}_b$ is completely analogous.

\begin{lemma}
\label{lemma-space-distance}
Let $v_1, v_2, \dots, v_m$ and $v'_1, v'_2, \dots, v'_{m'}$ form bases of subspaces $E$ and $E'$ of $V$, respectively.
Define matrices $F$, $G$, $H$ as follows,
$$
  F = \left( a(v_i,v'_j) \right)_{m\times m'},\quad  
  G = \left( a(v_i,v_j) \right)_{m\times m},\quad 
  H = \left( a(v'_i,v'_j) \right)_{m'\times m'}.
$$
Then
$$
\hat{\varepsilon}_a^2(E,E')=\lambda_{max}(F^TG^{-1}F, H)=\lambda_{max}(FH^{-1}F^T, G)\:,
$$
where $\lambda_{max}(A,B)$ denotes the maximum eigenvalue of eigen-problem $Ax=\lambda Bx$.

Further, suppose $\|F^TF\|_2\le \eta_F$, $\|I-G\|_2 \le \eta_G$, $\|I-H\|_2 \le \eta_H$. If $\eta_G, \eta_H <1$, then
$$
\hat{\varepsilon}_a^2(E,E') \le \frac{\eta_F}{(1-\eta_G)(1-\eta_H)}\:.
$$
\end{lemma}
\begin{proof}
Expand $v\in E$ and $v'\in E'$ as $v=\sum_{i=1}^m c_i v_i$ and $v'=\sum_{j=1}^{m'} c'_j v'_j$
and consider vectors $\bc\in\R^m$ and $\bcp\in\R^{m'}$ of coefficients $c_i$ and $c'_j$, respectively.
Then   
$$
  a(v,v') = \bc^T F \bcp, \quad
  \norm{v}_a^2 = \bc^T G \bc, 
  \quad\text{and}\quad
  \norm{v'}_a^2 = (\bcp)^T H \bcp\:.
$$
Thus, definition \eqref{eq:defzeta} gives
$$
  \hat{\varepsilon}_a(E,E') 
  = \max_{\bc^T G \bc = 1} \max_{(\bcp)^T H \bcp = 1} \bc^T F \bcp
  = \max_{\bc^T G \bc = 1} \max_{|\bctp| = 1} \bc^T F L^{-T} \bctp
    = \max_{\bc^T G \bc = 1} | \bc^T F L^{-T} |,
$$
where $\bctp = L^T \bcp$, $H = LL^T$ is the Cholesky decomposition of matrix $H$, and $|\cdot|$ stands for the Euclidean norm.
Consequently,
$$
\hat{\varepsilon}_a^2(E,E')  
= \max_{0 \neq \bc \in \R^{m}} \frac{\bc^T F L^{-T} L^{-1} F^T \bc}{\bc^T G \bc}
= \lambda_{max}(FH^{-1}F^T, G).
$$
Expression $\hat{\varepsilon}_a^2(E,E') = \lambda_{max}(F^TG^{-1}F, H)$ can be proved analogously.

To prove the upper bound on $\hat{\varepsilon}_a$, we use Cholesky decomposition $G=QQ^T$. 
Noticing that $\|A\|_2=\|A^T\|_2=\sqrt{\|A^TA\|_2}$ holds for a general matrix $A$, we have
$$
\lambda_{max}(FH^{-1}F^T, G) = \|Q^{-1}FH^{-1}F^T Q^{-T} \|_2
\le \|G^{-1}\|_2  \|H^{-1}\|_2  \|F^{T}F\|_2.
$$

Finally, we estimate $\|G^{-1}\|_2$ and $\|H^{-1}\|_2$.
If $\eta_G <1$, then
$$
\|G^{-1}\|_2 = \frac{1}{\lambda_{min} (G)}= \frac{1}{1-\lambda_{max} (I-G)} 
\le \frac{1}{1 - \|I-G\|_2} 
\le \frac{1}{1 - \eta_G}. 
$$
With the same argument for $H^{-1}$, we easily draw the conclusion.
\end{proof}

\begin{remark}
This lemma is used for spaces $\hE_k$ and $\hE_{k'}$ of approximate eigenfunctions with $k \neq k'$, therefore
matrices $F$, $G$, and $H$ are available and $\lambda_{max}(F^TG^{-1}F, H)$ as well as $\lambda_{max}(FH^{-1}F^T, G)$ can be computed. Alternatively, guaranteed estimates $\eta_F$, $\eta_H$, and $\eta_G$ can be obtained by the Gershgorin circle theorem. 
These estimates are expected to be good for $\hat{\varepsilon}_a(\hE_k,\hE_{k'})$,
because if the approximate eigenfunctions in $\hE_k$ and the ones in $\hE_{k'}$ are appropriately orthonormalized, then
$F^T F \approx 0$, $G\approx I_{m}$, and $H\approx I_{m'}$.
\end{remark}

The following theorem provides the desired estimates of the directed distances $\delta_a(E_K, \hE_K)$ and $\delta_b(E_K, \hE_K)$  defined in \eqref{eq:Delta}.

\begin{theorem}
\label{th:mainenergy}
Let the above specified partition of the spectrum into $K$ clusters be arbitrary.
Let $\hat u_i \in V$ for $i=1,2,\dots,N$ with $N < d = \operatorname{dim} V$
be such that $\operatorname{dim}\hE_k = N_k - n_k + 1$ for all $k=1,2,\dots,K$.
Let $\lambda_n < \rho \leq \lambda_{N+1}$.
Then
\begin{align}
  \label{eq:Deltaest}
  \delta_a^2(E_K,\hE_K)
  &\leq  {\frac{\rho (\hat\lambda^{(K)}_N -  \lambda_n)  + \lambda_n \hat\lambda^{(K)}_N \theta^{(K)}_a}{\hat\lambda^{(K)}_N(\rho - \lambda_n)}}
\quad\text{and}
\\
  \label{eq:Deltaest_b}
  \delta_b^2(E_K,\hE_K)
  &\leq {\frac{\hat\lambda^{(K)}_N - \lambda_n + \theta^{(K)}_b}{\rho - \lambda_n}} ,
\end{align}
where 
\begin{gather*}
  \hat\lambda^{(K)}_N = \max_{\hat v \in \hE_K} \frac{\norm{\hat v}_a^2}{\norm{\gamma \hat v}_b^2},
\quad  
  \theta^{(K)}_a = \sum_{k=1}^{K-1} \frac{\rho - \lambda_{n_k}}{ \lambda_{n_k} } \left[ \hat{\varepsilon}_a(\hE_k,\hE_K) + \delta_a(E_k,\hE_k) \right]^2,
\quad\text{and}  
\\
  \theta^{(K)}_b = \sum_{k=1}^{K-1} \left(\rho - \lambda_{n_k}\right) \left[ \hat{\varepsilon}_b(\hE_k,\hE_K) + \delta_b(E_k,\hE_k) \right]^2.
\end{gather*}
\end{theorem}

\begin{proof}
Following \eqref{eq:Parseval}, the proof is based on the equality 
\begin{equation}
\rho \norm{\gamma \hat u}_b^2 - \norm{ \hat u}_a^2 + \|\hat{u}_0\|_a^2
  = \sum_{i=1}^d \frac{\rho - \lambda_i}{\lambda_i} a( \hat u,u_i) ^2 
\end{equation}
for an arbitrary and fixed $\hat{u} \in \hE_K$, where 
$\hat{u}_0$ is the component of $\hat{u}$ in $\operatorname{Ker}(T \gamma)$ and $d$ can be infinity. 
Recalling that $n=n_K$ and $N = N_K$, using definition
\begin{equation}
  \vartheta(\hat u)
  = \sum_{i=1}^{n-1} \frac{\rho - \lambda_i}{\lambda_i} a( \hat u, u_i )^2, 
\end{equation}
and inequalities $\|\hat{u}_0\|_a^2 \ge 0$ and $\rho \leq \lambda_{N+1}$, we easily derive bound
$$
\rho \norm{\gamma \hat u}_b^2 - \norm{ \hat u}_a^2 - 
  \vartheta(\hat u) 
  \le \sum_{i=n}^N \frac{\rho - \lambda_i}{\lambda_i} a(\hat u, u_i)^2.
$$
Introducing orthogonal projectors $P_k: V \rightarrow E_k$, $k=1,2,\dots,K$,  
with respect to the inner product $a(\cdot , \cdot)$,
recalling their properties
\begin{equation}
\label{eq:tildeu}
\norm{P_k \hat u}_a^2 = \sum_{i=n_k}^{N_k} a( \hat{u}, u_i)^2,
\quad
\|\gamma P_k \hat u\|_b^2 = \sum_{i=n_k}^{N_k} \frac{a(\hat{u}, u_i)^2}{\lambda_i},
\end{equation}
and using $\lambda_n \leq \lambda_i$ for $i=n,\dots,N$, we immediately obtain estimates
\begin{equation}
\label{eq:expr_by_projection}
\rho \norm{\gamma \hat u}_b^2 - \norm{ \hat u}_a^2 - \vartheta(\hat u) 
  \le (\rho - \lambda_n) \norm{\gamma P_K \hat u}_b^2
  \le \frac{\rho - \lambda_n}{\lambda_n} \norm{P_K \hat u}_a^2.
\end{equation}
Similarly, the fact that $\lambda_{n_k} \leq \lambda_i$ for $i=n_k, \dots, N_k$ implies the following estimate on $\vartheta(\hat u)$:
\begin{equation}
\label{eq:theta_by_projection}
\vartheta(\hat{u}) 
\leq \sum_{k=1}^{K-1} \left(\rho-\lambda_{n_k}\right) \norm{\gamma P_k \hat u}_b^2
\leq \sum_{k=1}^{K-1} \frac{\rho - \lambda_{n_k}}{\lambda_{n_k}} \norm{P_k \hat u}_a^2.
\end{equation}

Now, we bound $\norm{P_k\hat{u}}_a$ and $\norm{\gamma P_k\hat{u}}_b$ for $k=1,2,\dots,K-1$.
Introducing $z_k = P_k \hat u / \norm{ P_k \hat u}_a \in E_k$ and
the $a$-orthogonal projector $\hP_k : V \rightarrow \hE_k$, 
definition \eqref{eq:defzeta} and relation \eqref{eq:Delta2} imply
$$
|a( \hat{u}, \hP_k z_k ) | \le \hat{\varepsilon}_a(\hE_k,\hE_K) \norm{\hat{u}}_a \norm{\hP_k z_k}_a
\quad\text{and}\quad
\norm{z_k - \hP_k z_k }_a \le \delta_a(E_k,\hE_k).
$$
Since $P_k \hat u = a( \hat{u}, z_k) z_k$ and
$\norm{ \hP_k z_k}_a \le \norm{ z_k}_a = 1$, these estimates provide the bound
\begin{multline*}
  \norm{ P_k \hat u}_a 
  = | a( \hat{u}, z_k ) | 
  \le | a( \hat{u}, \hP_k z_k ) | + | a( \hat{u}, z_k - \hP_k z_k ) | 
 \\ 
  \le \left[ \hat{\varepsilon}_a(\hE_k,\hE_K) + \delta_a(E_k,\hE_k) \right] \norm{\hat{u}}_a.
\end{multline*}

Analogous steps yield the bound
\begin{equation*}
  \norm{ \gamma P_k \hat u}_b 
  \le \left[ \hat{\varepsilon}_b(\hE_k,\hE_K) + \delta_b(E_k,\hE_k) \right] \norm{\gamma \hat{u}}_b.
\end{equation*}
Consequently, estimates \eqref{eq:theta_by_projection} provide 
\begin{equation}
\label{eq:theta_K_est}
\vartheta(\hat{u}) \le  
\theta^{(K)}_a \norm{\hat{u}}_a^2
\quad\text{and}\quad
\vartheta(\hat{u}) \le \theta^{(K)}_b \norm{\gamma \hat{u}}_b^2.
\end{equation}
The desired lower bounds on $\norm{P_K\hat{u}}_a$ and $\norm{\gamma P_K\hat{u}}_b$ follow from
\eqref{eq:expr_by_projection} and \eqref{eq:theta_K_est}:
\begin{align*}
  \| P_K \hat{u} \|_a^2  &\ge \lambda_n \frac{\rho \norm{\gamma \hat u}_b^2 - \norm{\hat u}_a^2 - \theta^{(K)}_a \norm{\hat{u}}_a^2 }{ \rho - \lambda_{n} },
\\ 
  \| \gamma P_K \hat{u} \|_b^2  &\ge \frac{\rho \norm{\gamma \hat u}_b^2 - \norm{\hat u}_a^2 - \theta^{(K)}_b \norm{\gamma\hat{u}}_b^2 }{ \rho - \lambda_{n} }.
\end{align*}
The final step is to apply these two estimates in expressions \eqref{eq:Delta3} for the distances $\delta_a^2(\hE_K,E_K)$ and $\delta_b^2(\hE_K,E_K)$.
Indeed,
$$
  \delta_a^2(\hE_K,E_K) = 1 - \min_{\substack{\hat{u} \in \hE_K\\ \norm{\hat{u}}_a = 1}} \norm{P_K \hat{u}}_a^2
  \leq
  1 + \max_{\substack{\hat{u} \in \hE_K\\ \norm{\hat{u}}_a = 1}} \lambda_n \frac{-\rho \norm{\gamma \hat u}_b^2 + \norm{\hat u}_a^2 + \theta^{(K)}_a \norm{\hat{u}}_a^{2}}{ \rho - \lambda_{n} }
$$
and the statement \eqref{eq:Deltaest} follows by elementary manipulations utilizing the definition of $\hat\lambda^{(K)}_N$. Statement \eqref{eq:Deltaest_b} follows analogously.
Note that $\delta_a(E_K,\hE_K) = \delta_a(\hE_K,E_K)$ and $\delta_b(E_K,\hE_K) = \delta_b(\hE_K,E_K)$, 
because $\operatorname{dim} E_K = \operatorname{dim} \hE_K = N - n + 1$.
\end{proof}

\begin{remark}\label{rem:birkhoff-etl-result}
Theorem~\ref{th:mainenergy} is
a direct and nontrivial generalization of \cite[Corollary~1]{BirBooSwaWen:1966}. 
Indeed, if all eigenvalues $\lambda_i$ are simple and well separated (forming clusters of size one),  $E_i = \linspan \{u_i\}$ and  $\hE_i = \linspan \{\hat u_i\}$ stand for the corresponding one-dimensional exact and approximate eigenspaces, $\hat\lambda_i = \norm{\hat u_i}_a^2/\norm{\hat u_i}_b^2$, and the corresponding approximate eigenfunctions $\hat{u}_i$, $i=1,2,\dots,k$, are mutually orthogonal,
then bounds \eqref{eq:Deltaest} and \eqref{eq:Deltaest_b} become
$$
  \delta_a^2(E_k,\hE_k)
  \leq  {\frac{1}{\rho - \lambda_k}}\left(\rho \frac{\hat\lambda_k -  \lambda_k}{\hat\lambda_k} + \lambda_k
  \sum_{i=1}^{k-1} \frac{\rho - \lambda_{i}}{ \lambda_{i} } \delta_a^2(E_i,\hE_i)
  \right)
$$
and
$$
  \label{eq:Deltaest_b_simple}
  \delta_b^2(E_k,\hE_k)
  \leq \frac{1}{{\rho - \lambda_k}} \left(\hat\lambda_k - \lambda_k + \sum_{i=1}^{k-1} \left(\rho - \lambda_{i}\right) \delta_b^2(E_i,\hE_i) \right),
$$
where the estimate for $\delta_b$ coincides with the statement in \cite[Corollary 1]{BirBooSwaWen:1966}.
\end{remark}

\begin{remark}
\label{remark:necessity-of-term-in-main-theorem}
In this remark, we show the necessity of terms  
$\rho - \lambda_n$, 
$\hat\lambda_N^{(K)} - \lambda_n$ and $\theta_a^{(K)}$ or $\theta_b^{(K)}$
appearing in bounds \eqref{eq:Deltaest} and \eqref{eq:Deltaest_b}.

\begin{itemize}
    \item [(i)] The difference $\rho - \lambda_n$ is determined by the spectral gap between the last cluster and the following eigenvalues. To have a sharp bound of the approximation error, the gap should not be too small.
\item [(ii)] 
Quantity $\hat\lambda_N^{(K)} - \lambda_n$ corresponds to the width of the last cluster. 
Note that in case of non-degenerated cluster, i.e., $\lambda_{N} > \lambda_{n}$, 
this difference causes bounds \eqref{eq:Deltaest} and \eqref{eq:Deltaest_b}
not to converge to zero even if $\hat\lambda_N^{(K)}\to \lambda_N$.
In fact, under the assumption that only the two-side bound of the whole cluster is available, the distribution of exact eigenvalues is uncertain and so are the corresponding exact eigenvectors. 
To illustrate this uncertainty in eigenvectors, let us consider a matrix eigenvalue problem $Ax=\lambda x$ with $A$ selected as
\begin{equation}
\label{eq:matrix-eig-example}
    A=\left(\begin{array}{cccc}1-\beta &0 &0 &0 \\ 0 & 1+\beta-\alpha & 0 & \sqrt{\alpha^2+2\alpha} \\ 0 & 0 & 1+\beta & 0 \\ 0& \sqrt{\alpha^2+2\alpha}&0& 3+\beta-\alpha \end{array} \right), \quad \beta > \alpha>0, \beta\approx 0\:.
\end{equation}

Eigenvalues of $A$ are given by 
$$
\lambda_1 = 1-\beta,\quad \lambda_2=1 +\beta-2\alpha,\quad \lambda_3=1+\beta,\quad \lambda_4 = 3+\beta,
$$
and the corresponding 
the eigenvectors are
$$
u_1=e_1,\quad 
u_2= \left( \begin{array}{c}0 \\ \frac{-\sqrt{\alpha+2}}{\sqrt{2\alpha+2}} \\ 0\\ \frac{\sqrt{\alpha}}{\sqrt{2\alpha+2}}\end{array} \right),\quad 
u_3=e_3,\quad 
u_4=\left( \begin{array}{c} 0  \\  \frac{\sqrt{\alpha}}{\sqrt{2\alpha+2}} \\ 0\\ \frac{\sqrt{\alpha+2}}{\sqrt{2\alpha+2}}  \end{array} \right).
$$
Let $a(u,v):=u^TAv$, $b(u,v):=u^Tv$ for column vector $u,v \in \R^4$.
Considering the approximate eigenvectors 
$\widehat{u}_i=e_i$ for $i=1,2,3$ and spaces 
$E_1 = \linspan\{u_1, u_2,u_3\}$,
$\widehat{E}_1 = \linspan\{\widehat{u}_1,\widehat{u}_2,\widehat{u}_3\}$, 
we compute the distance $\delta_b^2(E_1, \widehat{E}_1)$ exactly and apply bound \eqref{eq:Deltaest_b} with $\rho=\lambda_4$ to obtain 
$$
\delta_b^2(E_1, \widehat{E}_1) 
\left( = \frac{\alpha}{2+2\alpha}\right)\le 
\frac{2\beta}{2+2\beta}.
$$
Clearly, the quality of approximate vectors measured by the distance $\delta_b$ depends on the position of $\lambda_2$ (determined by $\alpha$). For such a system with uncertain $\alpha$, one cannot expect a sharper estimation for $\delta_b$.

If further information about eigenvalues inside a tight cluster is known, we can  (theoretically) split the cluster into smaller clusters consisting 
of a single or a multiple eigenvalue and for these clusters bounds \eqref{eq:Deltaest} 
and \eqref{eq:Deltaest_b} do converge. 
Bounds \eqref{eq:Deltaest} and \eqref{eq:Deltaest_b} are naturally computed iteratively starting from the first cluster. Accuracy of this procedure is illustrated on numerical examples  in Section~\ref{se:numex}.
\item [(iii)] 

Values of $\theta_a^{(K)}$ and $\theta_b^{(K)}$ measure errors in all previous clusters. 
Notice that quantity $\theta_a^{(K)}$ (and similarly $\theta_b^{(K)}$) depends on $\delta_a(E_k,\hE_k)$ 
and $\hat\varepsilon_a(\hE_k,\hE_K)$ for $k=1,2,\dots,K-1$.
The distance $\delta_a(E_k,\hE_k)$ accounts for errors in spaces of eigenfunctions of previous clusters
and $\hat\varepsilon_a(\hE_k,\hE_K)$ for possible non-orthogonality of approximate eigenfunctions.

Note that the dependence of bounds \eqref{eq:Deltaest} and \eqref{eq:Deltaest_b} on errors in previous clusters is necessary, because these bounds utilize only the information about eigenvalue gaps and Rayleigh quotients.
Indeed, consider the one-dimensional Dirichlet eigenvalue problem: 
\begin{equation}
    \label{eq:eig-lambda-1d}
-u^{''}=\lambda u \mbox{ in } I=(0,\pi);~~ u(0)=u( \pi )=0.
\end{equation}
 and approximate eigenfunctions 
$$
\hat{u}_1(x) = \sin x, \quad
\hat{u}_2(x) = (1-t) \sin x  + t \sin (3x), \text{ with }  t = {\sqrt{3}}/({\sqrt{3}+\sqrt{5}}). 
$$

Clearly, $\hat{u}_2$ is a bad approximation of the exact eigenfunction $u_2(x) = \sin (2x)$.
Taking naturally $E_2 = \linspan\{u_2\}$, $\hE_2 = \linspan\{\hat u_2\}$, distances $\delta_a^2(E_2,\hE_2)$ and $\delta_b^2(E_2,\hE_2)$ are relatively large.
On the other hand, the corresponding Rayleigh quotient gives the exact eigenvalue, namely $R(\hat{u}_2) = \hat\lambda_2^{(2)} = \lambda_2 = 4$. 
The difference $\hat\lambda_2^{(2)} - \lambda_2$ in \eqref{eq:Deltaest} and \eqref{eq:Deltaest_b} is zero and values of these upper bounds are determined by terms $\theta_a^{(2)}$ and $\theta_b^{(2)}$ only, reflecting the non-orthogonality of $\hat{u}_1$ and $\hat{u}_2$.

Terms $\theta_a^{(2)}$ and $\theta_b^{(2)}$ could be avoided if, for example, guaranteed estimates of the residual are employed. This approach will be introduced in \S\ref{subsec:extension-of-davis-kahan} as Algorithm II.

\end{itemize}
\end{remark}

The following theorem addresses the question of efficiency of bounds presented in Theorem~\ref{th:mainenergy}. To formulate it, we denote the right-hand sides of estimates \eqref{eq:Deltaest} and \eqref{eq:Deltaest_b} by $\overline{\delta}_a(E_K, \widehat{E}_K)$ and $\overline{\delta}_b(E_K, \widehat{E}_K)$, respectively.
\begin{theorem}
\label{th:efficiency-analysis}
If eigenspaces $\widehat{E}_k$ for $k=1,2,\dots,K$ are mutually orthogonal 
then 
there exists a generic constant $C > 0$ determined by the leading ${N_K}$ exact and approximate eigenvalues and independent of the width of clusters such that
\begin{equation}\label{eq:eff}
\max\left\{ \overline{\delta}^2_a(E_K, \widehat{E}_K),
\overline{\delta}^2_b(E_K, \widehat{E}_K) \right\}
\le C \sum_{k=1}^K \left({\delta}^2_a(E_k, \widehat{E}_k)  + \lambda_{N_k} - \lambda_{n_k}\right).
\end{equation}
\end{theorem}

\begin{proof}
Take $\hat{u} \in \widehat{E}_K$ such that $\|\hat{u}\|^2_a=\hat{\lambda}^{(K)}_N \|\gamma \hat{u}\|_b^2$. 
Let $u=P_K \hat{u}\in E_K$ be the best approximation to $\hat{u}$ under the norm $\|\cdot\|_a$. By Lemma~\ref{le:delta}, 
$\|\hat{u}-u\|_a \le \delta_a(E_K, \widehat{E}_K)$. 
Suppose
$u=\sum_{i=n}^N c_i u_i$ and take $\tilde{\lambda}\in [\lambda_n,\lambda_N]$ such that 
$\tilde{\lambda} \sum_{i=n}^N c_i^2 = \sum_{i=n}^N c_i^2\lambda_i$.
By expanding the squared terms $\|\hat{u}-u\|_a^2$ and $\|\gamma(\hat{u}-u)\|_b^2$, we have
$$
\|\hat{u}-u\|_a^2 - \tilde{\lambda} \|\gamma(\hat{u}-u)\|_b^2
= 
2 b\left(\sum_{i=n}^N c_i (\tilde{\lambda} - {\lambda_i})  u_i, \hat{u}\right)
+(\hat{\lambda}_N^{(K)} - \tilde{\lambda}) \|\gamma\hat{u}\|_b\:.
$$
Noticing that $\|\gamma(\hat{u}-u)\|_b^2 \le 1/\lambda_1 \|\hat{u}-u\|_a^2$ and 
$\lambda_n \le \tilde{\lambda} \le \lambda_N$, the equation above tells
\begin{equation}
\label{eq:local_est_for_eigen_approximate}
|\hat{\lambda}_N^{(K)} - \lambda_N| \le C (\|\hat{u}-u\|_a^2 +
\lambda_N -\lambda_n) \:.
\end{equation}
Therefore, for the term $(\hat{\lambda}_N^{(K)}-\lambda_n)$ in (\ref{eq:Deltaest}) and (\ref{eq:Deltaest_b}), we have
$$
\hat{\lambda}_N^{(K)}-\lambda_n = \hat{\lambda}_N^{(K)}-\lambda_N + \lambda_N-\lambda_n \le 
C (\delta_a^2(E_K,\widehat{E}_K) +  \lambda_{N} - \lambda_{n})  \:.
$$
The mutual orthogonality of $\widehat{E}_k$ for $k=1,2,\dots,K$ implies 
$\hat{\epsilon}_a(E_k, E_K)=\hat{\epsilon}_b(E_k, E_K)=0$. 
Consequently, terms $\theta_a^{(K)}$ and $\theta_b^{(K)}$ can be estimated recursively, leading to the inequality \eqref{eq:eff}.
\end{proof}

\medskip

As pointed out in (ii) of Remark \ref{remark:necessity-of-term-in-main-theorem}, the term $\lambda_{N_k} -\lambda_{n_k} $ 
is necessary for clusters of non-zero width.
If all clusters are of zero size and errors in previous clusters are comparable to the error in the cluster of interest, i.e., 
$$
  \delta_a^2(E_k,\widehat E_k) \leq C \delta_a^2(E_K,\widehat E_K)
  \quad\forall k=1,2,\dots,K-1,
$$ 
then inequality \eqref{eq:eff} immediately yields the following efficiency result
$$
\max\left\{ \overline{\delta}^2_a(E_K, \widehat{E}_K),
\overline{\delta}^2_b(E_K, \widehat{E}_K) \right\}
\le C {\delta}^2_a(E_K, \widehat{E}_K).
$$

\begin{remark}
In Section~\ref{se:numex}, we apply the general bounds \eqref{eq:Deltaest} and \eqref{eq:Deltaest_b} in the particular context of the finite element method. This method, specifically, computes the best approximation ($a$-orthogonal projection) of the exact solution of the underlying boundary value problem in the finite element space. As the numerical results in Section~\ref{se:numex} reveal, the bound \eqref{eq:Deltaest} is optimal in the sense that it has the same rate of convergence as $\delta_a(\hE_K, E_K)$. 
However, the bound \eqref{eq:Deltaest_b} has a lower rate of convergence than $\delta_b(\hE_K, E_K)$.

The Aubin--Nitsche technique together with explicit {\em a priori} error estimates \cite{LiuOis2013} will provide fully computable and optimal bounds under $b$-norm. However, since the 
application of the Aubin--Nitsche technique utilizes  special properties about FEM approximation, such an approach is beyond the problem setting of this paper and will be addressed in a subsequent publication; partial results of this approach can be found in \cite{liu-Vejchodsky-arxiv-v1}.
\end{remark}

\medskip

\section{Sharp $a$-norm estimates based on the $b$-norm bounds}
\label{se:a_by_b}
This section provides an estimate of the distance $\delta_a(E,\hE)$ by the distance $\delta_b(E,\hE)$ between the space of exact eigenfunctions $E$ and the space of approximate eigenfunctions $\hE$.
The idea is motivated by the following well known formula (see e.g. \cite[page 55]{Boffi:2010})
\begin{equation}
\label{eq:fundamental-relation-between-eigenfunction-and-eigenvalue}
 \| u_i - \hat{u}_i\|_a^2 = \lambda_i \|\gamma(u_i - \hat{u}_i) \|_b^2 - (\lambda_i - \hat\lambda_i) \| \gamma \hat{u}_i \|_b^2
\end{equation}
for the exact eigenpair $\lambda_i > 0$, $u_i \in V$, arbitrary approximate eigenfunction $\hat{u}_i \in V$ and approximate eigenvalue $\hat\lambda_i = \| \hat{u}_i \|_a^2 / \| \gamma \hat{u}_i \|_b^2$.
In the context of the finite element method for the Laplace eigenvalue problem, this identity essentially says that the error 
$\|u_i - \hat{u}_i\|_a$ in the energy norm is dominated by the error of the approximate eigenvalue itself, because 
the error $\|\gamma (u_i - \hat{u}_i)\|_b$ in the $L^2$ norm has a higher order of convergence.

The following estimate is theoretically independent of the partition of eigenvalues into clusters,
but its natural usage is to bound $\delta_a(E_k, \hE_k)$ by $\delta_b(E_k, \hE_k)$, 
where $k$ is the index of a cluster as it is introduced at the end of Section~\ref{se:eigenproblem}.

\begin{theorem}
\label{th:a_by_b}
Let $u_n, \dots, u_N$ be the exact eigenfunctions of \eqref{eq:eigp} and $0 < n \leq N$ the corresponding indices.
Let $\hat u_n, \dots, \hat u_N \in V$ be linearly independent.
Let $E = \linspan\{u_n, \dots, u_N\}$ and $\hE = \linspan\{\hat u_n, \dots, \hat u_N\}$.
Then
\begin{equation}
  \label{eq:a_by_b}
  \delta_a^2(E,\hE) \leq {2 - 2 \lambda_n \left( \frac{1 - \delta_b^2(E,\hE) }{\lambda_N \hat\lambda_N} \right)^{1/2}}, 
\end{equation}
where $\lambda_n$ and $\lambda_N$ are exact eigenvalues corresponding to $u_n$ and $u_N$ and 
$$
\hat\lambda_N = \max_{\hat v \in \hE} \frac{\norm{\hat v}_a^2}{\norm{\gamma \hat v}_b^2}.
$$
\end{theorem}
\begin{proof}
Consider the linear mapping $\tau: E \to E$ defined by
$$
\tau (u) = \sum_{i=n}^N c_i\lambda_i u_i,\quad
\text{where } u=\sum_{i=n}^N c_i u_i.
$$
Since $\lambda_i > 0$ for all $i=n,\dots,N$, 
$\tau$ is a bijection. Given arbitrary $u\in E$ and $\hat{u}\in\hE$, we clearly have
$$
a(u, \hat{u}) 
  = \sum_{i=n}^N c_i a(u_i, \hat{u}) 
  = \sum_{i=n}^N c_i \lambda_i b(\gamma u_i , \gamma \hat{u})
  = b(\gamma \tau(u), \gamma \hat{u}).
$$
This enables us to estimate the distance \eqref{eq:Delta} between $E$ and $\hE$ as follows
\begin{multline}
\label{eq:Deltaest_a}
\delta_a^2(E, \hE)  
  = \max_{ \substack{u \in E \\ \norm{u}_a = 1} } \min_{\hat{u} \in \hE} \| u - \hat{u} \|_a^2 
  \le \max_{ \substack{u \in E \\ \norm{u}_a = 1} } \min_{ \substack{ \hat{u} \in \hE \\ \norm{\hat{u}}_a = 1} } \|u - \hat{u} \|_a^2 
\\
  = \max_{\substack{u \in E \\ \norm{u}_a = 1}} \min_{\substack{ \hat{u} \in \hE \\ \norm{\hat{u}}_a = 1}} \left[ 2 - 2 b(\gamma \tau(u), \gamma \hat{u}) \right]
 \leq 2 - 2 \lambda_n \min_{\substack{u \in E \\ \norm{u}_a = 1}} \max_{\substack{ \hat{u} \in \hE \\ \norm{\hat{u}}_a = 1}} b\left(\frac{\gamma \tau(u)}{\|\tau(u)\|_a}, \gamma \hat{u} \right),
\end{multline}
where the last inequality follows from the fact that
$$
  \norm{\tau(u)}_a^2 = \sum_{i=n}^N \lambda_i^2 c_i^2 \geq \lambda_n^2 \sum_{i=n}^N c_i^2 = \lambda_n^2 \norm{u}_a^2 = \lambda_n^2
  \quad\forall u \in E,\ \norm{u}_a=1.
$$
Since $\tau$ is a bijection, it is easy to show that
$$
    \left\{ \frac{\tau(u)}{\norm{\tau(u)}_a} : u\in E,\ \norm{u}_a = 1 \right\}
  = \left\{ u \in E : \norm{u}_a = 1 \right\}.
$$
This equality together with bounds $\norm{u}_a^2 \leq \lambda_N \norm{\gamma u}_b^2$ for all $u\in E$ and $\norm{\hat{u}}_a^2 \leq \hat\lambda_N \norm{\gamma \hat{u}}_b^2$ for all $\hat{u} \in \hE$ imply
\begin{multline}
\label{eq:minmaxestb}
   \min_{\substack{u \in E \\ \norm{u}_a = 1}} \max_{\substack{ \hat{u} \in \hE \\ \norm{\hat{u}}_a = 1}} b\left(\frac{\gamma \tau(u)}{\norm{\tau(u)}_a}, \gamma \hat{u}\right)
  = \min_{\substack{u \in E \\ \norm{u}_a = 1}} \max_{\substack{ \hat{u} \in \hE \\ \norm{\hat{u}}_a = 1}} b(\gamma u, \gamma \hat{u})
\\  
  = \min_{\substack{u \in E \\ u\neq0}} \max_{\substack{\hat{u} \in \hE \\ \hat{u}\neq0}} b\left(\frac{\gamma u}{\norm{u}_a},  \frac{\gamma \hat{u}}{\norm{\hat{u}}_a}\right)
  = \min_{\substack{u \in E \\ \norm{\gamma u}_b=1}} \max_{\substack{\hat{u} \in \hE \\ \norm{\gamma \hat{u}}_b=1}} b\left(\frac{\gamma u}{\norm{u}_a},  \frac{\gamma \hat{u}}{\norm{\hat{u}}_a}\right)
\\  
  \geq \frac{1}{\left(\lambda_N \hat\lambda_N\right)^{1/2}} \min_{\substack{u \in E \\ \norm{\gamma u}_b=1}} \max_{\substack{\hat{u} \in \hE \\ \norm{\gamma \hat{u}}_b=1}} b(\gamma u, \gamma \hat{u})
  = \left( \frac{1 - \delta_b^2(E,\hE) }{\lambda_N \hat\lambda_N} \right)^{1/2},
\end{multline}
where we note that $\max_{\hat{u} \in \hE,\ \norm{\gamma \hat{u}}_b=1} b(\gamma u, \gamma \hat{u})$ is non-negative and the last equality follows from \eqref{eq:Delta3} using characterization 
$$
  \norm{\gamma \hP u}_b = \max_{\substack{\hat{u} \in \hE \\ \norm{\gamma \hat{u}}_b=1}} b(\gamma u, \gamma \hat u).
$$
The proof is finished by substituting \eqref{eq:minmaxestb} to \eqref{eq:Deltaest_a}.
\end{proof}

\begin{remark}
The estimate \eqref{eq:a_by_b} will provide a sharper bound on $\delta_a(E,\widehat{E})$
if the eigenvector approximation has a higher convergence rate under $\|\cdot\|_b$ norm compared with $\|\cdot\|_a$ norm, which usually happens for the finite element approximation of elliptic eigenvalue problems.

For a multiple eigenvalue, we have $\lambda_n = \lambda_N$ and the right-hand side of \eqref{eq:a_by_b} is close to 
$$ 
2 \frac{\hat\lambda_N^{1/2} - \lambda_N^{1/2}}{\hat\lambda_N^{1/2}} 
+ \frac{\lambda_N^{1/2}}{\hat\lambda_N^{1/2}} \delta_b^2(E,\hE)
= 2 \frac{\hat\lambda_N - \lambda_N}{\hat\lambda_N^{1/2}\left(\hat\lambda_N^{1/2} + \lambda_N^{1/2}\right)}
+ \frac{\lambda_N^{1/2}}{\hat\lambda_N^{1/2}} \delta_b^2(E,\hE),
$$
where we use the approximation $(1-\delta_b^2)^{1/2} \approx 1-\delta_b^2/2$.
This shows that for example in the context of the finite element method for the Laplace eigenvalue problem,
the bound \eqref{eq:a_by_b} has the optimal rate of convergence. 
Numerical experiments indicate that the bound \eqref{eq:a_by_b} combined with the sub-optimal estimate \eqref{eq:Deltaest_b} can still provide a sharper bound on $\delta_a$ than \eqref{eq:Deltaest}; see Section~\ref{se:numex}.
\end{remark}

Define the following quantity to measure the distance between $E$ and $\hE$:
\begin{equation}
\label{eq:def-tilde-delta-as-the-final-select}
\tilde{\delta}(E, \hE )
:=\max_{ \substack{\hat{u} \in \hE \\ \norm{\gamma \hat{u}}_b = 1} }  \min_{\substack{ {u} \in E }} \| u- \hat{u}\|_a.
\end{equation}
This quantity can be easily bounded by $\delta_a$ as follows
\begin{equation}
\label{eq:relation_of_tilde_delta_and_delta_a}
\sqrt{\hat{\lambda}}_n \delta_a (E, \hE ) \le 
\tilde{\delta}(E, \hE ) \le \sqrt{\hat{\lambda}}_N \delta_a (E, \hE ).
\end{equation}
The following lemma relates $\tilde{\delta}$ and $\delta_b$.
\begin{lemma}
\label{lem:relation_of_tilde_delta_and_delta_b}
By using $\delta_b(E,\hE)$, we have the following estimate for 
$\tilde{\delta}(E, \hE)$,
\begin{equation}
\label{eq:est-tilde-delta-by-delta-b}
\tilde{\delta}^2(E, \hE)
\le 
\lambda_N + \hat{\lambda}_N - 2 \lambda_n \sqrt{1-\delta_b^2(E,\hE)}.
\end{equation}

\end{lemma}
\begin{proof}
We utilize the same mapping $\tau: E \to E$ as in the proof of Theorem \ref{th:a_by_b} to obtain the estimate for $\tilde{\delta}$.
Note that for any $ u \in E$, $\norm{\gamma u}_b=1$, we have $\norm{\gamma \tau(u)}_b \ge  \lambda_n$ and consequently
\begin{align*}
\tilde{\delta}^2(E, \hE)   
&  \le \max_{ \substack{\hat{u} \in \hE \\ \norm{\gamma \hat{u}}_b = 1} } \min_{ \substack{ {u} \in E \\ \norm{\gamma u}_b = 1} } \left[ \lambda_N +\hat{\lambda}_N - 2 b(\gamma \tau(u), \gamma \hat{u}) \right] \notag
\\
& \leq \lambda_N +\hat{\lambda}_N - 2 \lambda_n \min_{ \substack{\hat{u} \in \hE \\ \norm{\gamma\hat{u}}_b = 1} } \max_{ \substack{ {u} \in E \\ \norm{\gamma u}_b = 1} } b\left(\frac{\gamma \tau(u)}{\|\gamma \tau(u)\|_b}, \gamma \hat{u} \right)
	\\
	& =\lambda_N +\hat{\lambda}_N - 2\lambda_n \min_{ \substack{\hat{u} \in \hE \\ \norm{\gamma \hat{u}}_b = 1} } \max_{ \substack{ {u} \in E \\ \norm{\gamma u}_b = 1} } b(\gamma u, \gamma \hat{u})\\
	& =  \lambda_N +\hat{\lambda}_N - 2 \lambda_n \left( {1 - \delta_b^2(E,\hE) } \right)^{1/2}~.
\end{align*}
\end{proof}

\begin{remark}
If approximate eigenvalues converge to the exact ones, i.e. if $|\lambda_{k,h}-\lambda_k| \to 0$, then the right-hand sides of estimates \eqref{eq:a_by_b} and \eqref{eq:est-tilde-delta-by-delta-b} do not converge to zero, due to the existence of the term corresponding to the cluster width $|\lambda_N -\lambda_n|$. The necessity of this term is confirmed by the discussion of the eigenvalue problem for matrix \eqref{eq:matrix-eig-example}, where $\delta_a \to \sqrt{2\beta +\beta^2}$  and $\delta_b \to 0$ when $\alpha \to 0$.
Note that $\beta$ determines the cluster width for the problem in \eqref{eq:matrix-eig-example}.
Despite the inefficiency caused by the cluster width, 
the advantage of utilizing (\ref{eq:est-tilde-delta-by-delta-b}) is that for clusters with narrow  width, an easy-to-obtain sharp bound on
$\delta_b$ will greatly improve the precision for $\hat{\delta}$ and thus $\delta_a$.
\end{remark}

\section{Algorithm II: Estimation based on the residual error}
\label{se:est-residual-error}

Estimates in Theorem~\ref{th:mainenergy} of \S \ref{se:ebound} only utilize the information about eigenvalues distribution and Rayleigh quotients for approximate eigenfunctions. As discussed in Remark \ref{remark:necessity-of-term-in-main-theorem},
such estimates
have the advantage to provide guaranteed {error bound} for problems with uncertainty information about the objective eigenvalue problem. 
As a result, the bounds in (\ref{eq:Deltaest}) and (\ref{eq:Deltaest_b}) depend on the quality of previous clusters and the width of the cluster of interest.

In this section, we utilize the residual error estimation of approximate eigenfunctions to obtain  sharper {bounds} of the distances between the exact and the approximate eigenspaces.
For this purpose, we introduce the approach based on the Davis--Kahan's method \cite{davis1970rotation}.
Note that the original Davis--Kahan's method only deals with strongly formulated differential operators and we extend it to weakly formulated eigenvalue problems.
It is worth to point out that, compared to very cheap estimates given in Theorem~\ref{th:mainenergy}, the approach utilizing the residual will require more computational effort to compute or bound its suitable norm.
The complexity of this effort varies and depends on the problem.

The Davis--Kahan's $\sin\theta$ theorem in \cite{davis1970rotation}  bounds the error of an approximate eigenspace by a norm of the residual and a spectral gap. The bound requires to apply the operator to an approximate eigenfunction to estimate the norm of the residual. This can be easily done for strongly defined operators and sufficiently smooth approximations. However, for weakly defined operators and non-smooth approximations the original Davis--Kahan's $\sin\theta$ theorem cannot be directly used. Therefore, we extend the Davis--Kahan's $\sin\theta$ theorem to weakly defined eigenvalue problems in the form \eqref{eq-main-eig-form}.

For the reader's convenience, we quote Davis--Kahan's result \cite{davis1970rotation} using the notation of this paper. 

\begin{theorem}[{Davis-Kahan's $\sin \theta$ theorem}]
\label{thm:davis-kahan-origin}
Let $U$ be a Hilbert space with inner product $\langle \cdot, \cdot \rangle$. Let $\|\cdot\|$ be a unitary invariant norm of $U$. 
Let $A:U\to U$ be a self-adjoint operator and $E$ be the subspace spanned by the eigenvectors of $A$ in the cluster of interest. That is,
$$
  {E}=\linspan\{u_{n},\dots, u_{N}\}, \quad
  A u_i = \lambda_i u_i \quad {\text{for } i=n,\dots,N}.
$$
Let $A_h$ be a perturbation of the operator $A$, which has an invariant subspace $\widehat{E}$ such that 
$\widehat{E}=\linspan\{\hat{u}_n,\dots, \hat{u}_{N}\}$
and 
$$
  A_h \hat{u}_i = \hat\lambda_i \hat{u}_{i}\quad{\text{for } i=n,\dots, N}.
$$
Let $P:U\to E$ and $\widehat{P}:U\to \widehat{E}$ be the two projection operators with respect to $\langle \cdot, \cdot \rangle$. 
Suppose the eigenvalues  of $A$ corresponding to $E^\perp$  are excluded from 
$(\hat\lambda_n-\sigma, \hat\lambda_N+\sigma)$.
Then
$$
\|(I-P)\widehat{P}\| \le \frac{\|(I-P)(A-A_h)\widehat{P}\|}{\sigma} \le \frac{\|(A-A_h)\widehat{P}\|}{\sigma}~.
$$
\end{theorem}

Since the proof in the original paper is divided into several parts for more general discussion, here we provide a reformulation of the proof in \cite{davis1970rotation} in a concise way. 

\begin{proof}
Define 
$c=(\hat\lambda_n +\hat\lambda_N )/2$, $r=(\hat\lambda_N - \hat\lambda_n)/2 $.
Since $AP=PA$, we have
\begin{align}
(I-P)(A-A_h)\widehat{P} &=(I-P)A\widehat{P} - (I-P) A_h \widehat{P} \notag\\
&= (A-cI)(I-P)\widehat{P} - (I-P)(A_h-cI)\widehat{P} \:. \label{eq:proof-davis-kahan-1}
\end{align}
Since the eigenvalues of $A-cI$  over $E^\perp$ are excluded from $(-r-\sigma, r+\sigma)$,  we have 
\begin{equation}
\label{eq:proof-davis-kahan-2}
\|(I-P)\widehat{P} \|= \| (A-cI)^{-1}|_{E^\perp} (A-cI)(I-P)\widehat{P}  \| \le \frac{1}{r+\sigma} \| (A-cI)(I-P)\widehat{P}  \|~.
\end{equation}

Also, since the eigenvalue of $A_h-cI$  over $\widehat{E}$ belongs to $[-r,r]$, we have
\begin{equation}
\label{eq:proof-davis-kahan-3}
\|(I-P)(A_h-cI)\widehat{P}\| = \|(I-P)\widehat{P} (A_h-cI)\widehat{P}\|  \le  r \|(I-P)\widehat{P}\|~.
\end{equation}

Thus, we draw the conclusion from (\ref{eq:proof-davis-kahan-1}), (\ref{eq:proof-davis-kahan-2}), (\ref{eq:proof-davis-kahan-3}).
\end{proof}

\begin{remark}
If the operator $A$ is strongly defined as a differential operator, then for any given smooth test function $u_h$ from $U$, we can easily evaluate $A u_h$.  However, if the problem is weakly formulated, one cannot evaluate $Au_h$ directly. In the succeeding subsection, we will introduce the extension of Davis--Kahan's result to weakly formulated eigenvalue problems. 
\end{remark}

\subsection{Extension of Davis--Kahan's $\sin\theta$ theorem to weakly formulated problems }
\label{subsec:extension-of-davis-kahan}

Let us recall the weakly formulated eigenvalue problem \eqref{eq-main-eig-form} and the composition $T\gamma : V \to V$ of the compact operator $\gamma$ and the solution operator $T$.  
We set $U = \mbox{Ker}(T\gamma)^\perp = \mbox{Ker}(\gamma)^\perp(\subset V)$ and notice that the operator $T\gamma$ is invertible on $U$. 
Thus, the operator $A=\left(T\gamma|_U\right)^{-1} : U \to U$ is well defined and satisfies
$$
  a(u, v) = b(\gamma Au, \gamma v) \quad \forall u,v \in U.
$$
Since $T\gamma$ is self-adjoint, the operator $A$ is self-adjoint as well and we may apply the Davis--Kahan's method for it.

As above, we consider an eigenvalue cluster of interest and denote the space of exact eigenfunctions corresponding to eigenvalues in the cluster by $E = \linspan\{u_n, \dots, u_N\}$. 
These eigenfunctions are approximated by the one from the finite dimensional space 
$\hE = \linspan\{\hat u_n, \dots, \hat u_N \}$.
Note that $\hE$ need not be necessarily a subspace of $U$. 
This setting enables us to define an approximate operator $A_h:\hE \to \hE$ by the identity
$$
  a(\hat u, \hat v) = b(\gamma A_h \hat u, \gamma \hat v) \quad\forall \hat u, \hat v \in \hE.
$$
Thus, the discrete eigenvalue problem to find $\hat\lambda_i$ and $\hat u_i \in \hE \setminus\{0\}$ such that
$$
  a(\hat u_i, \hat v) = \hat\lambda_i b(\gamma \hat u_i, \gamma \hat v) \quad \forall \hat v \in \hE,
$$
is equivalent to the eigenvalue problem 
$$
  A_h \hat u_i = \hat\lambda_i \hat u_i \quad {\text{for } i=n,\dots, N}.
$$

Recall that $P$ and $\widehat{P}$ are projections that map $V$ to $E$ and $\hE$ with respect to the inner product $a(\cdot, \cdot)$, respectively.
To measure the residual error of $(A-A_h)$, let us introduce the quantity $\epsilon$ by
\begin{equation}
\label{def:epsilon-for-residual}
\epsilon:=\frac{1}{\sigma} \sup_{u \in U, \|\gamma u\|_b=1 } \|(I-P)(A-A_h)\widehat{P} u\|_{-1}~.	
\end{equation}
Here, the norm $\|\cdot\|_{-1}$ is defined as 
\begin{equation}
    \label{def:inverse_norm}
    \|u\|_{-1} := \max_{v\in U}\frac{b(\gamma u,\gamma v)}{\|v\|_a} \quad \text{for } u \in U~.
\end{equation}
Note that
Definitions 3.4, 3.5 from 
\cite[\S 3.2]{CanDusMadStaVoh2019} and Definition 2.1 from \cite{CanDusMadStaVoh2017}
use a similar quantity to $\epsilon$ to derive the estimator.

\begin{lemma}
\label{lem:davis-kahan-proof-norm-scaling}
Given a normed space $B$ with a norm $\|\cdot\|_B$, let $L:B \to B$ be a linear mapping such that $\|L\|_B \le r$. {Let $F:B\to \mathbb{R}$ be an absolutely homogeneous functional (i.e., $F(\alpha v) = |\alpha| F(v)$ for all $\alpha \in \R$ and $v \in B$).} Then 
$$
\sup_{v\in B, \|v\|_B=1 } F(L v) \le r \sup_{v\in B, \|v\|_B=1 } F(v)~.
$$
\end{lemma}
\begin{proof} The proof is easily done by noting the following {inequality},
\begin{align*}
  \sup_{\substack{v\in B \\ \|v\|_B=1 }} F(L v) &
  \leq \sup_{v\in B} \frac{F(L v)}{\norm{Lv}_B} \frac{\norm{Lv}_B}{\norm{v}_B} 
  \leq \sup_{y\in B} \frac{F(y)}{\norm{y}_B} \norm{L}_B 
  = \norm{L}_B \sup_{\substack{v\in B\\ \|v\|_B=1}} F(v),	
\end{align*}
where the first inequality follows from the absolute homogeneity of $F$ and the second one from the fact that $\operatorname{Range}(L) \subset B$ and from the definition of the operator norm. The final equality uses again the absolute homogeneity of $F$.
\end{proof}

\medskip

We formulate the extension of the Davis--Kahan's $\sin\theta$ theorem
below in Theorem~\ref{thm:extend-of-davis-kahan-for-delta-b}.

\begin{theorem}
\label{thm:extend-of-davis-kahan-for-delta-b}
Suppose the eigenvalues of $A$ corresponding to $E^\perp$ are excluded from
$(\hat\lambda_n-\sigma, \hat\lambda_N+\sigma)$
i.e., assume that there exists $\sigma > 0$ such that 
$$
\lambda_{n-1}\leq \hat\lambda_n-\sigma < \hat\lambda_N+\sigma \leq \lambda_{N+1}.
$$ 
Then
\begin{equation}
\label{eq:relation-delta_b_and_tilde_delta}
\delta_b(E,\hE) \le \epsilon \cdot \tilde{\delta} (E,\hE).	
\end{equation}

\end{theorem}

\begin{proof}
Define 
$c:=(\hat\lambda_n +\hat\lambda_N )/2$, $r:=(\hat\lambda_N - \hat\lambda_n)/2$.
The proof is based on the triangle inequality
\begin{equation}
\label{eq:davis-kahan-extension-1}    
\|(I-P)(A-A_h)\widehat{P} u\|_{-1} \ge \|(A-cI)(I-P)\widehat{P}u\|_{-1} - \|(I-P)(A_h-cI)\widehat{P}u\|_{-1},
\end{equation}
where $u\in U$.
Note for any $\widehat{u}$ in $\hE$, 
$(A_h -cI)\widehat{u}\in \hE$ and 
$\|\gamma (A_h -cI)\widehat{u}\|_b \le r \|\gamma \widehat{u}\|_b$.
Lemma~\ref{lem:davis-kahan-proof-norm-scaling} leads to the following inequality,
\begin{equation}
\label{eq:davis-kahan-extension-2}    
\sup_{u \in U, \|\gamma u\|_b=1 } \|(I-P)(A_h-cI)\widehat{P}u\|_{-1} \le r \sup_{u \in U, \|\gamma u\|_b=1 } \|(I-P)\widehat{P} u\|_{-1}~.
\end{equation}
Note that for any $v \in E^\perp$, 
$\|\gamma (A-cI) v \|_b \ge (r+\sigma) \|\gamma v\|_b$. 
For any $u \in E^\perp$, by taking $v:=(A-cI)|_{E^\perp}^{-1}u \in E^\perp$, we have
$$
\|\gamma u\|_b=\|\gamma (A-cI)  \cdot 
  {(A-cI)|_{E^\perp}^{-1} u}
\|_b
\ge (r+\sigma) \|\gamma (A-cI)|_{E^\perp}^{-1}u\|_b.
$$
That is,
\begin{equation}
    \label{eq:local_inv_A_minus_I_est}
\|\gamma (A-cI)|_{E^\perp}^{-1}u\|_b \le \frac{1}{(r+\sigma)} \|\gamma u\|_b~.
\end{equation}

Noticing that $(I-P)\widehat{P}  = (A-cI)|_{E^\perp}^{-1} (A-cI) (I-P)\widehat{P}$,  Lemma \ref{lem:davis-kahan-proof-norm-scaling} along with \eqref{eq:local_inv_A_minus_I_est} implies
\begin{equation}
\label{eq:davis-kahan-extension-3}
\sup_{u \in U, \|\gamma u\|_b=1 } \|(A-cI)(I-P)\widehat{P}u\|_{-1} \ge (r+\sigma)  \sup_{u \in U, \|\gamma u\|_b=1 } \|(I-P)\widehat{P} u\|_{-1}.
\end{equation}

Inequalities (\ref{eq:davis-kahan-extension-1}),  (\ref{eq:davis-kahan-extension-2}), and  (\ref{eq:davis-kahan-extension-3}) lead to the following result.
$$
\sup_{u \in U, \|\gamma u\|_b=1 } \|(I-P)\widehat{P} u\|_{-1} \le \frac{1}{\sigma} \sup_{u \in U, \|\gamma u\|_b=1 } \|(I-P)(A-A_h)\widehat{P} u\|_{-1} =\epsilon~.
$$
Next,  we apply the above estimation to bound $\delta_b(E,\hE)$. Note that
\begin{align*}
    \|\gamma (I-P)\widehat{P} u)\|_b & = \sup_{v\in U, \|\gamma v\|b =1 }  b(\gamma (I-P)\widehat{P} u), \gamma (I-P)\widehat{P} v)\\ 
    &=\sup_{v\in U, \|\gamma v\|b =1}   \|(I-P)\widehat{P} u\|_{-1} \cdot  \|(I-P)\widehat{P} v\|_a.
\end{align*}

Therefore, we obtain the following estimation for $\delta_b(E,\hE)$:
\begin{align*}
~\delta_b(E,\hE)  &= \sup_{u\in U, \|\gamma u\|_b=1}  \|\gamma (I-P)\widehat{P} u\|_b  \\
& \le \sup_{u\in U, \|\gamma u\|_b=1}  \|(I-P)\widehat{P} u\|_{-1} ~\cdot ~ \sup_{v\in U, \|\gamma v\|_b=1}\|(I-P)\widehat{P} v\|_a \\
& \le \epsilon \cdot \tilde{\delta} (E,\hE)~.
\end{align*}

\end{proof}

Applying relation \eqref{eq:est-tilde-delta-by-delta-b} between $\delta_b$ and $\tilde{\delta}$ to the statement of Theorem~\ref{thm:extend-of-davis-kahan-for-delta-b},
we obtain the following {bound on} $\delta_b(E,\hE)$.
\begin{equation}
    \label{eq:estimation-of-b-norm-by-inequality}
{\delta_b}(E,\hE)^2  \le \epsilon^2 \cdot \left(\lambda_N+\hat{\lambda}_N -2\lambda_n \sqrt{1-\delta_b(E,\hE)^2}\right)~.
\end{equation}

Inequality (\ref{eq:estimation-of-b-norm-by-inequality}) can be reformulated as a quadratic inequality for $t:=\delta_b^2$,
$$
g(t):=t^2 +\left(\alpha_1^2 -2\alpha_2 \right)t  + \alpha_2^2 - \alpha_1^2 \ge 0,\quad  t \le \alpha_2,
$$
with
$\alpha_1 = 2\epsilon^2 \lambda_n$ and 
$\alpha_2 = \epsilon^2 \cdot (\lambda_N+\hat{\lambda}_N)$.
By solving {this quadratic inequality,} 
we {derive} the following theorem.
\begin{theorem}
\label{thm:extension-davis-kahan-for-delta-b}
Let the two solution of $g(t)=0$ be $t_1, t_2$. In case $\alpha_1 < \alpha_2$ and $0< t_1 \le \alpha_2 < t_2$, then we have
$\delta_b^2 \le t_1$. That is, 
\begin{equation}
    \label{eq:estimate-delta-b-by-residual-info}
\delta_b \le \left\{\frac{2\alpha_2-\alpha_1^2 - \alpha_1  \sqrt{\alpha_1^2+4\alpha_2-4} }{2}
\right\}^{1/2}.
\end{equation}
\end{theorem}
{Estimates} for $\delta_a$ and $\tilde{\delta}$ {are} available through the relation \eqref{eq:a_by_b} and \eqref{eq:est-tilde-delta-by-delta-b}.

\begin{remark}
Note that for small value of  $\delta_b$, we have $\sqrt{1-\delta_b^2}\approx (1-\delta_b^2/2)$. Hence 
$$
\delta_b \lessapprox \epsilon \left\{ \frac{ \lambda_N +\hat{\lambda}_N- 2\lambda_n}{1-\epsilon^2 \lambda_n}\right\}^{1/2} = \epsilon \cdot \left(O( \mbox{cluster width})+O(|\hat{\lambda}_N-\lambda_N|)\right)^{1/2}~.
$$	
The quality of the estimator for $\delta_b$ depends on the residual error, the width of the specified cluster and the quality of the eigenvalue estimation.
As we will see in next sub-section, the quantity $\epsilon$ does not depend on the width of the cluster. In solving practical problems,
if the width of the specified cluster is large, 
 one can try {to split} the cluster to reduce the cluster width.

 In the case of solving the Laplace eigenvalue problem by the conforming linear finite element method, for which it is expected that 
 $O(|\hat{\lambda}_N-\lambda_N|)=O(h^2)$, 
 $\epsilon = O(h)$ and $\delta_b = O(h^2)$. Our proposed estimation coincides with this fact, if the cluster width is zero or small enough compared {to} $O(|\hat{\lambda}_N-\lambda_N|)$.

\end{remark}

\begin{remark}
The following identity for any $u_h \in \hE$ and $v\in U$, will be useful for computing a tight bound on $\epsilon$:	
\begin{align}
  b(\gamma (I-P)(A-A_h)u_h, \gamma v)    & = b( \gamma (A-A_h)u_h, \gamma (I-P)v)  
  \nonumber \\
     &  = a(u_h, (I-P)v) - b(\gamma A_h u_h , \gamma (I-P)v)~.
  \label{eq:residual_expression}
\end{align}
\end{remark}

\subsection{Weakly formulated residual error estimation}

In this subsection, we consider a concrete eigenvalue problem of the Laplacian and utilize the Prager--Synge technique \cite{PraSyn:1947}  to compute a tight bound on the residual error term $\epsilon$. The technique relies on a flux reconstruction, which we explain for the case of the Laplace eigenvalue problem with homogeneous Dirichlet boundary condition. 
The Prager--Synge technique has wide applications in both {\em a posteriori} and {\em a priori} error estimation; see for example \cite{Liu-2022} for a recent successful application in validating the solution to the Navier--Stokes equation in 3D domains. 

\medskip

Let $V = H^1_0(\Omega)$, $W = L^2(\Omega)$, $\gamma$ be the identity, $a(u,v) = (\nabla u, \nabla v)$, and $b(\gamma u, \gamma v) = (u,v)$, where the parenthesis denote the $L^2(\Omega)$ scalar product. Further, let $\|\cdot\|_0$ stands for the $L^2(\Omega)$ norm.

\begin{lemma}
\label{le:eqfluxest}
Let $u_h \in H^1_0(\Omega)$ and $p \in \Hdiv$ be arbitrary. Then
$$
  \|(I-P)(A-A_h)u_h\|_{-1} \leq \|\nabla u_h-p\|_0 + \frac{1}{\sqrt{\lambda_1}} \| \ddiv p +A_h u_h\|_0.
$$
\end{lemma}
\begin{proof}
Using the particular forms of $a$ and $b$ for the Laplace eigenvalue problem, the negative norm in the definition of $\epsilon$ can be expressed with the aid of \eqref{eq:residual_expression} as
\begin{equation}
\label{eq:negnormest}
\|(I-P)(A-A_h)u_h\|_{-1} = \max_{v \in V} \frac{(\nabla u_h, \nabla (I-P)v) - (A_h u_h, (I-P)v)}{\norm{\nabla v}_0}~.
\end{equation}
Denoting $\varphi = (I-P)v$, choosing arbitrary $p\in\Hdiv$, applying the divergence theorem, Cauchy--Schwarz inequality, and the Friedrichs inequality $\norm{\varphi}_0 \leq \norm{\nabla\varphi}_0 / \sqrt{\lambda_1}$, we obtain
\begin{align*}
  (\nabla u_h, \nabla \varphi) - (A_h u_h, \varphi) 
&  = (\nabla u_h - p, \nabla \varphi) - (\ddiv p + A_h u_h, \varphi) 
\\
&  \leq \norm{\nabla u_h - p}_0 \norm{\nabla\varphi}_0 + \norm{\ddiv p + A_h u_h}_0 \norm{\varphi}_0 
\\  
&  \leq \left( \norm{\nabla u_h - p}_0  + \frac{1}{\sqrt{\lambda_1}}\norm{\ddiv p + A_h u_h}_0\right) \norm{\nabla\varphi}_0.  
\end{align*}
Inserting this estimate to \eqref{eq:negnormest} together with inequality $\norm{\nabla\varphi}_0 \leq \norm{\nabla v}_0$ finishes the proof.
\end{proof}
\begin{remark}
The residual error $\|(I-P)(A-A_h)u_h\|_{-1}$ is the same as the one in \cite[Definition 2.1]{CanDusMadStaVoh2017}, up to the unimportant factor $I-P$. However, the assumptions and error bounds for approximate eigenfunctions obtained in \cite{CanDusMadStaVoh2017} are different from our approach. We assume high-precision bounds for eigenvalues in the cluster of interest, which can be obtained by applying, for example, the Lehmann--Goerisch theorem with the FEM over refined meshes; while \cite{CanDusMadStaVoh2017} considers bounds for eigenvalues and eigenvectors simultaneously. Numerical results in \S\ref{se:numex} 
compare the accuracy of these two approaches and show that our approach is capable of high precision.
\end{remark}

\begin{lemma}
\label{le:epsmaxmin}
Let $W_h$ be a finite dimensional subspace of $\Hdiv$ and $\hE$ the space of approximate eigenvalues corresponding to the cluster of interest, then
\begin{equation}
\label{eq:estimate-for-espilon}
  \epsilon \leq \frac{1}{\sigma} \max_{\substack{u_h \in \hE \\ \|u_h\|_0=1 }} ~  \min_{p_h \in W_h } 
  \left( \|\nabla u_h-p_h\|_0 + \frac{1}{\sqrt{\lambda_1}} \| \ddiv p_h +A_h u_h\|_0 \right).
\end{equation}
Further, if $A_h \hE \subset \ddiv W_h$, we have
\begin{equation}
  \label{eq:estimate-for-espilon-simple}
\epsilon \leq \frac{1}{\sigma} \max_{\substack{u_h \in \hE \\ \|u_h\|_0=1 }} ~
\min_{ \substack{p_h \in W_h,\\ \ddiv p_h +A_h u_h=0}} \|\nabla u_h -p_h \|_0. 
\end{equation}
\end{lemma}
\begin{proof}
Since $\hP  V = \hE$ and using Lemma~\ref{le:eqfluxest}, we have
\begin{multline*}
  \sup_{u \in V, \|u\|_0=1 } \|(I-P)(A-A_h)\hP u\|_{-1} 
  = \max_{u_h \in \hE, \|u_h\|_0=1 } \|(I-P)(A-A_h)u_h\|_{-1} 
\\
\le \max_{u_h \in \hE, \|u_h\|_0=1 } \min_{p_h \in W_h } \|\nabla u_h-p_h\|_0 
    + \frac{1}{\sqrt{\lambda_1}} \| \ddiv p_h + A_h u_h\|_0.
\end{multline*}
In case that $A_h u_h \in \ddiv W_h$, we can thus force $\ddiv p_h + A_h u_h=0$ in seeking the minimizer $p_h$.
Inserting this bound to the definition of $\epsilon$ finishes the proof.
\end{proof}

Note that in the case of the finite element method the approximate eigenvectors $u_h$ forming $\hE$ are piecewise polynomial. A natural choice for $W_h$ in Lemma~\ref{le:epsmaxmin} is a suitable Raviart--Thomas space of higher order such that we can always find $p_h \in W_h$ satisfying $\ddiv p_h +A_h u_h=0$. There are  alternative methods to obtain this flux; see, e.g., \cite{CanDusMadStaVoh2017,Rep-2008-book}.
Assuming this choice, we consider linear operator $R_h: \hE \to W_h$ mapping $u_h$ to $p_h = R_h u_h$, which solves the following minimization problem:
$$
\min_{p_h \in W_h,\ \ddiv p_h +A_h u_h=0} \|\nabla u_h -p_h \|_0. 
$$
With the operator $R_h$, the evaluation of the residual error term $\epsilon$ reduces to solving a  maximization problem
$$
\epsilon \le \frac{1}{\sigma} \max_{u_h \in \hE }\frac{ \|\nabla u_h -R_h u_h\|_0}{\|u_h\|_0}.
$$
The right-hand side of the above inequality corresponds to an eigenvalue problem for matrices with dimension $\operatorname{dim}(\hE)$.

\begin{remark}
The bound in Theorem \ref{thm:extension-davis-kahan-for-delta-b}
utilizes relations \eqref{eq:a_by_b} and \eqref{eq:est-tilde-delta-by-delta-b} to obtain estimates for $\delta_a$ and $\tilde{\delta}$. Notice that a large cluster width $\lambda_N - \lambda_n$ will cause a drop of precision of this bound. 
Therefore, the following subsection derives bounds on $\delta_a$ and $\tilde{\delta}$ independent from the cluster width. We achieve it by introducing another application of Davis--Kahan's theorem and estimating $\delta_a$ and $\tilde{\delta}$ directly. 
\end{remark}

\subsection{{Direct estimation of $\delta_a$: another application of  Davis--Kahan's theorem}} 
Corresponding to the operator $T\gamma$ defined in (\ref{eq:def-T-operator}), let us introduce a local discrete operator $\widehat{T}{\gamma}:\hE\to\hE$. For given $\hat{u} \in \hE$, the image $\widehat{T}{\gamma} \hat{u}$ is uniquely determined as the solution to problem
$$
a(\widehat{T}{\gamma} \hat{u}, \hat{v}) = b(\gamma \hat{u}, \gamma \hat{v}) \quad \forall \hat{v} \in \hE.
$$
Note that $\widehat{T}{\gamma}\hat{u}_i=\hat{\mu}_i \hat{u}_i$ for $i=n,\dots, N$, where 
$\hat{u}_i$ and $\hat{\mu}_i$ ($\hat{\mu}_i\ge \hat{\mu}_{i+1}$) are the discrete eigenpairs over $\hE$, i.e., 
$b(\gamma \hat{u}_i, \gamma \hat{v})=\hat{\mu}_i a(\hat{u}_i, \hat{v})$ for all $\hat{v} \in \hE$.
Suppose the eigenvalues of $T$ corresponding to $E^\perp$ are excluded from 
$(\hat{\mu}_N-\sigma, \hat{\mu}_n+\sigma)$, then from Davis--Kahan's theorem we have
$$
\delta_a(E,\hE) \le \frac{1}{\sigma} \|(I-P)(T{\gamma}-\widehat{T}{\gamma})\widehat{P}\|_a
$$
Note the following variational equation holds  for any $\hat{u}\in \hE$ {and $v\in V$},
$$
a((T {\gamma}-\widehat{T} {\gamma})\hat{u}, v) = b({\gamma} \hat{u}, {\gamma} v) - a(\widehat{T}{\gamma} \hat{u}, v)~.
$$
Then, we have
\begin{align*}
\|(I-P)(T{\gamma}-\widehat{T}{\gamma})\widehat{P}\|_a  & \le \|(T{\gamma}-\widehat{T}{\gamma})\widehat{P}\|_a
=\sup_{u\in V, \|u\|_a=1} \|(T{\gamma}-\widehat{T}{\gamma})\widehat{P}u\|_a \\
&=\sup_{\hat{u}\in \hE, \|\hat{u}\|_a=1} \sup_{v\in V, \|v\|_a=1} b({\gamma}\hat{u}, {\gamma} v) - a(\widehat{T}{\gamma} \hat{u}, v)~.
\end{align*}
In case of the Laplace eigenvalue problem, suppose $\hE$ is the approximation by using a conforming FEM space and $W_h$ the FEM approximation to $H(\mbox{div})$, for example, the Raviart--Thomas space. Then,  with analogous argument as in Lemma \ref{le:epsmaxmin}, we have
\begin{equation}
\label{eq:estimate-for-a-norm-by-davis-kahan-direct}
\delta_a(E,\hE) \leq \frac{1}{\sigma} \max_{\substack{u_h \in \hE \\ \|\nabla u_h\|_0=1 }} ~  \min_{p_h \in W_h } 
  \left( \|\nabla \widehat{T}{\gamma} u_h-p_h\|_0 + \frac{1}{\sqrt{\lambda_1}} \| \ddiv p_h +  u_h\|_0 \right).
\end{equation}

Compared with the approach utilizing 
\eqref{eq:estimate-delta-b-by-residual-info} and \eqref{eq:a_by_b} for estimating $\delta_a$, 
the estimate \eqref{eq:estimate-for-a-norm-by-davis-kahan-direct} is independent from the cluster width $|\lambda_N-\lambda_n|$. Such an estimate is expected to deal with clusters of relative large width. 
{Note that, by further utilizing the relation in (\ref{eq:relation_of_tilde_delta_and_delta_a}) and (\ref{eq:relation-delta_b_and_tilde_delta}), 
the estimate \eqref{eq:estimate-for-a-norm-by-davis-kahan-direct} also applies to 
$\tilde{\delta}$ and thus for $\delta_{b}$. That is, 
$$
\delta_b \le \epsilon \, \tilde{\delta} 
\le \hat{\lambda}_N \, \epsilon \, \delta_a.
$$
This estimate of $\delta_b$ has the convergence rate as $O(h^2)$ for liner conforming FEM approximation if solution $u$ has the $H^2$ regularity.

\begin{remark}
\label{remark:super-convergence}
In practical computations with narrow cluster widths, 
the estimate  \eqref{eq:estimate-for-a-norm-by-davis-kahan-direct}
 gives relatively inaccurate bounds on $\delta_a$ in comparison with 
 Theorem 
 \ref{thm:extension-davis-kahan-for-delta-b}. 
This is because the estimate
\eqref{eq:estimate-for-a-norm-by-davis-kahan-direct}
introduces to $\delta_a$ an additional overestimating factor of convergence rate $O(h)$.
On the other hand,
\eqref{eq:estimate-for-espilon} 
introduces an overestimating factor of rate $O(h)$ for $\epsilon$ and 
the approach through Theorem 
 \ref{thm:extension-davis-kahan-for-delta-b}, which depends on 
 \eqref{eq:estimate-delta-b-by-residual-info} and \eqref{eq:a_by_b},  
finally leads to an overestimating factor for $\delta_a$ with the convergence rate $O(h^2)$. Note that the estimator for $\delta_a$ itself is still of convergence rate $O(h)$.
\end{remark}

\medskip

To summarize, we extended the Davis--Kahan theorem for weakly formulated eigenvalue problems and derived bounds \eqref{eq:estimate-delta-b-by-residual-info} for narrow clusters and \eqref{eq:estimate-for-a-norm-by-davis-kahan-direct} for wide clusters.
 These estimates are independent from the cluster index, its width and degree. The quality of these bounds is determined by the residual error and the gap between the concerned eigenvalue cluster and its neighbouring clusters.
{Thus, these estimates enjoy the same favourable properties as those in \cite{CanDusMadStaVoh2019}}.
As demonstrated by numerical examples in the following section, the estimates using the residual information provide sharper bound for $\delta_a$ and $\delta_b$ than the results in \cite{CanDusMadStaVoh2019}, especially for finite element approximations for eigenvalue problems of differential operators. 

\section{Numerical examples}
\label{se:numex}

This section numerically illustrates the accuracy of proposed bounds on the directed 
distances of spaces of exact and approximate eigenfunctions for matrix and Laplace eigenvalue problems.
Refer to Remark~\ref{re:examples} for how these problems fit to the considered abstract context.\footnote{See \url{https://ganjin.online/xfliu/EigenVecEstimation} for source codes and demonstrations of all presented examples.}

\subsection{The matrix eigenvalue problem}
\label{sec:num-matrix-eig}
First, we apply estimates of the directed distance between spaces of eigenvectors 
to the generalized matrix eigenvalue problem $Ax=\lambda Bx$, where $A, B \in \R^{961\times 961}$ are
the stiffness and mass matrices of the Laplacian in the unit square discretized by linear conforming finite element method on a uniform mesh with $h = \sqrt{2}/32$ (cf. Figure~\ref{fi:squaremesh}).
In agreement with the abstract eigenvalue problem \eqref{eq:eigp}, the $a$- and $b$- norms of a vector $v$ are given as 
$$
\|v\|_a = \sqrt{v^T A v}
\quad\text{and}\quad 
\|\gamma v\|_b = \sqrt{v^T B v}.
$$
The approximate eigenvectors are computed by the MATLAB command {\em eigs}. The error of these approximations is estimated by Theorem~\ref{th:mainenergy}. The needed two-sided bounds on eigenvalues are computed by using 
the interval arithmetic providing guaranteed 
enclosing intervals of exact eigenvalues.

The upper bounds on eigenvalues are easily obtained by the Rayleigh-Ritz method, 
while the lower bounds need more effort by applying Sylvester's law of inertia; 
see the detailed implementation in \cite{Behnke1991}. 

Table~\ref{tab:matrix_eig_pro_est}  lists the resulting two-sided bounds of eigenvalues and estimates of the directed distance of the corresponding subspaces of eigenvectors. 
Notice that since the approximate eigenvectors do no have better precision in the $b$-norm, the bound \eqref{eq:a_by_b} has no advantage over \eqref{eq:Deltaest_b}.

\begin{table}[ht]
    \centering
    \caption{Error estimates of eigenvalues and eigenvectors for the generalized matrix eigenvalue problem }
    \label{tab:matrix_eig_pro_est}
    \begin{tabular}{|c|l|l|l|l|l|}
    \hline
    \multirow{2}{*}{cluster}  & two-sided bounds         & eigenvalue           & bound \eqref{eq:Deltaest} & bound \eqref{eq:Deltaest_b} & bound \eqref{eq:a_by_b} \\
                              & on eigenvalues          & enclosure size       & on $\delta_a$       & on $\delta_b$ & on $\delta_a$     \\
    \hline
    \rule[-2mm]{0mm}{6mm}{}
    $\{\lambda_1\}$ & $19.78151183_{16}^{25}$ & $8.43\times10^{-10}$ & $4.10\times10^{-7}$ & $2.60\times10^{-7}$ & $6.23\times10^{-6}$\\ 
    \hline                              \rule[-2mm]{0mm}{6mm}{}
    $\{\lambda_{2},\lambda_{3}\}$ & $49.5769652_{58}^{60}$  & $1.00\times10^{-9}$  & $1.13\times10^{-6}$ & $6.12\times10^{-7}$ & $8.89\times10^{-7}$ \\
    \hline 
    \rule[-2mm]{0mm}{6mm}{}
    $\{\lambda_4\}$ & $79.63484906_{38}^{50}$ & $1.10\times10^{-9}$  & $5.33\times10^{-6}$ & $1.98\times10^{-6}$ & $1.20\times10^{-6}$ \\ 
    \hline                              \rule[-2mm]{0mm}{6mm}{}
    $\{\lambda_{5},\lambda_{6}\}$ & $99.50227701_{43}^{55}$ & $1.16\times10^{-9}$  & $1.06\times10^{-5}$ & $3.46\times10^{-6}$ & $2.04\times10^{-6}$ \\
    \hline
    \end{tabular}
\end{table}

\subsection{Comparison with the $\sin\theta$ theorem of Davis and Kahan} 
\label{sec:num-matrix-eig-standard}
Since the proposed bounds and Davis--Kahan's estimates can be easily compared in the case of the standard matrix eigenvalue problem, 
we consider the problem $Ax=\lambda x$, where $A$ is the same matrix as in Subsection~\ref{sec:num-matrix-eig}.
Its eigenvalues and eigenvectors can be easily computed with high accuracy. 
The leading $13$ eigenvalues of $A$ are naturally clustered as in Table \ref{table:cluster-matrix}.

\begin{table}[ht]
    \centering
    \caption{\label{table:cluster-matrix}The leading $8$ clusters for the eigenvalue problem $Ax=\lambda x$}
    \label{tab:standard_eig_pro_clusters}
    \begin{tabular}{|c|c|c|c|c|c|c|c|}
    \hline
    \rule[-1.5mm]{0mm}{5mm}{}
     $E_1$ &  $E_2$ &  $E_3$ &  $E_4$ &  $E_5$ & $E_6$ & $E_7$ & $E_8$  \\
     \hline
     \rule[-1.5mm]{0mm}{5mm}{}
    $\{\lambda_1\}$ &   $\{ \lambda_2, \lambda_3 \}$ & $\{\lambda_4\}$ & $\{\lambda_5,\lambda_6\}$ & $\{\lambda_7, \lambda_8\}$ & $\{\lambda_9, \lambda_{10}\}$ & $\{\lambda_{11}\}$ &  $\{\lambda_{12}, \lambda_{13}\}$ \\
     \hline
    \end{tabular}
\end{table}

To generate approximate eigenvectors $\hat u_i$, we perturb the exact eigenvectors $u_i$ by uniform random errors of magnitude $10^{-4}$ and apply the Gram--Schmidt process to satisfy the orthonormality requirement of the Davis--Kahan's method. 
Then, we compute approximate eigenvalues
by the Rayleigh quotient as $\hat{\lambda}_i = \hat{u}_i^T A \hat{u}_i / \hat{u}_i^T \hat{u}_i$
and observe that
differences $|\lambda_i - \hat{\lambda}_i|$ for the leading six eigenvalues are around $3\times10^{-6}$.

First, we estimate the distance $\delta_b(E_k, \hE_k)$ between the exact space $E_k$ and its approximation $\hE_k$ by the bound \eqref{eq:Deltaest_b}. It is evaluated by using two-sided bounds of eigenvalues of $A$ that are obtained independently by using the same method as for the generalized eigenvalue problem in Subsection~\ref{sec:num-matrix-eig}.
The $\sin \theta$ theorem of Davis and Kahan \cite{davis1970rotation} (i.e., Theorem \ref{thm:davis-kahan-origin})  provides the estimate
\begin{equation}
\label{eq:davis-kahan-est}
\delta_b(E_k, \hE_k) \le \frac{\|A \hat{\mathbf{x}} -\hat{\mathbf{x}} \hat{\Lambda} \|_{2}}{\delta_\mathrm{spec}},
\end{equation}
where $\hat{\mathbf{x}}=(\hat{u}_{n_k},\cdots, \hat{u}_{N_k})$ is the matrix with orthonormal column vectors that form  $\hE_k$; $\hat{\Lambda}$ is the diagonal matrix with diagonal elements $\hat{\lambda}_{n_k}, \cdots, \hat{\lambda}_{N_k}$; and
$\delta_\mathrm{spec}$ is the spectral gap between the cluster of interest and the rest of the spectrum. For example, in case $k=2$, we have $\hat{\mathbf{x}}=(\hat{u}_2, \hat{u}_3)$, $\hat{\Lambda}= \operatorname{diag}(\hat{\lambda}_2,\hat{\lambda}_3)$, and $\delta_\mathrm{spec}=\min(\lambda_4 - \lambda_3, \lambda_2 -\lambda_1)$. 
Table~\ref{tab:comparison-with-davis-kahan} compares the computed bounds with the exact values $\delta_b(E_k,\hE_k)$ for $k=1,2,\dots, 8$. 
The bound \eqref{eq:Deltaest_b} gives better results than the Davis--Kahan's estimate, except for the $7$th cluster $\{\lambda_{11}\}$. For the $8$th cluster, the Davis--Kahan's estimate 
exceeds one; note that the value of the directed distance cannot be greater than one.

\begin{table}[ht]
    \centering
    \caption{The directed distance $\delta_b(E_k, \hE_k)$, its bound \eqref{eq:Deltaest_b}, and Davis--Kahan's estimate \eqref{eq:davis-kahan-est} for clusters $k=1,2,\dots, 8$.}    \label{tab:comparison-with-davis-kahan}
    \begin{tabular}{|c|l|l|l|l|l|}    \hline     \rule[-2mm]{0mm}{6mm}{}
    Cluster  & $\{\lambda_1\} $   &  $\{\lambda_2, \lambda_3\} $ &   $\{\lambda_4\} $ &  $\{\lambda_5, \lambda_6\} $ \\
    \hline \rule[-2mm]{0mm}{6mm}{}
    $\delta_b(E_k, \hE_k)$  & 
    1.26 $\times 10^{-3}$ & 2.39 $\times 10^{-3}$ &  1.78 $\times 10^{-3}$ &  2.16 $\times 10^{-3}$ \\ 
    \hline \rule[-2mm]{0mm}{6mm}{}
    Bound \eqref{eq:Deltaest_b}  & 
     $1.08\times 10^{-2}$ & $ 2.12\times 10^{-2}$ & $4.36\times 10^{-2}$ & $7.24\times 10^{-2}$  \\
    \hline \rule[-2mm]{0mm}{6mm}{}
    Davis--Kahan \eqref{eq:davis-kahan-est} &
    $1.36\times 10^{-1}$ &$ 1.41\times 10^{-1} $ & $2.04\times 10^{-1}$ & $2.14\times 10^{-1}$  \\ 
    \hline 
    \end{tabular} \\
    \quad \\
    \begin{tabular}{|c|l|l|l|l|l|}
    \hline \rule[-2mm]{0mm}{6mm}{}
    Cluster  &
     $\{\lambda_7,\lambda_8\} $   &  $\{\lambda_9, \lambda_{10}\} $ &   $\{\lambda_{11}\} $ &  $\{ \lambda_{12}, \lambda_{13}\} $ \\
    \hline  \rule[-2mm]{0mm}{6mm}{}
    $\delta_b(E_k, \hE_k)$  & 
    2.33 $\times 10^{-3}$ & 2.34 $\times 10^{-3}$ &  1.83 $\times 10^{-3}$ &  2.40 $\times 10^{-3}$ \\ 
    \hline  \rule[-2mm]{0mm}{6mm}{}
    Bound \eqref{eq:Deltaest_b}  & 
     $1.26\times 10^{-1}$ & $ 3.72\times 10^{-1}$ & $5.64\times 10^{-1}$ & $8.53\times 10^{-1}$  \\
    \hline                                                \rule[-2mm]{0mm}{6mm}{}
    Davis--Kahan \eqref{eq:davis-kahan-est} &
    $2.18\times 10^{-1}$ &$ 3.86\times 10^{-1} $ & $3.99\times 10^{-1}$ & $18.4$  \\ 
    \hline 
    \end{tabular} 
\end{table}

\subsection{The unit square domain}

Consider the Laplace eigenvalue problem with homogeneous Dirichlet boundary conditions in the unit square $\Omega=(0,1)^2$: 
find eigenvalues $\lambda_i \in \R$ and corresponding eigenfunctions $u_i \neq 0$ such that
\begin{equation}
\label{eq:modpro}
  -\Delta u_i = \lambda_i u_i \quad\text{in }\Omega,
  \qquad
  u_i = 0 \quad\text{on }\partial\Omega.
\end{equation}
As we mention in Remark~\ref{re:examples}, the $a$- and $b$-norms correspond to the energy and $L^2(\Omega)$ norms, respectively.

The exact eigenpairs are known analytically to be
$$
\lambda_{ij} = (i^2+j^2)\pi^2,\quad u_{ij}=\sin(i\pi x) \sin(j\pi y), \quad i,j=1,2,3, \dots.
$$
These eigenvalues are either simple or multiple
and we clustered them according to the multiplicity.
The first four clusters are listed in Table~\ref{ta:square_domain_clusters}.
Since the exact eigenvalues are known, we use them to evaluate error bounds \eqref{eq:Deltaest}, \eqref{eq:Deltaest_b}, and \eqref{eq:a_by_b}. The quantity $\rho$ is chosen as $\rho = \lambda_{N+1}$. 
If the exact eigenvalues are not known, their two-sided bounds have to be employed.

This problem is discretized by the conforming finite element method of the first order. The finite element mesh $\cT_h$ is chosen as the uniform triangulation consisting of
isosceles right triangles; see the mesh with size $h=1/4$ in Figure~\ref{fi:squaremesh}.

\begin{figure}
\begin{center}
    \includegraphics[scale=0.4]{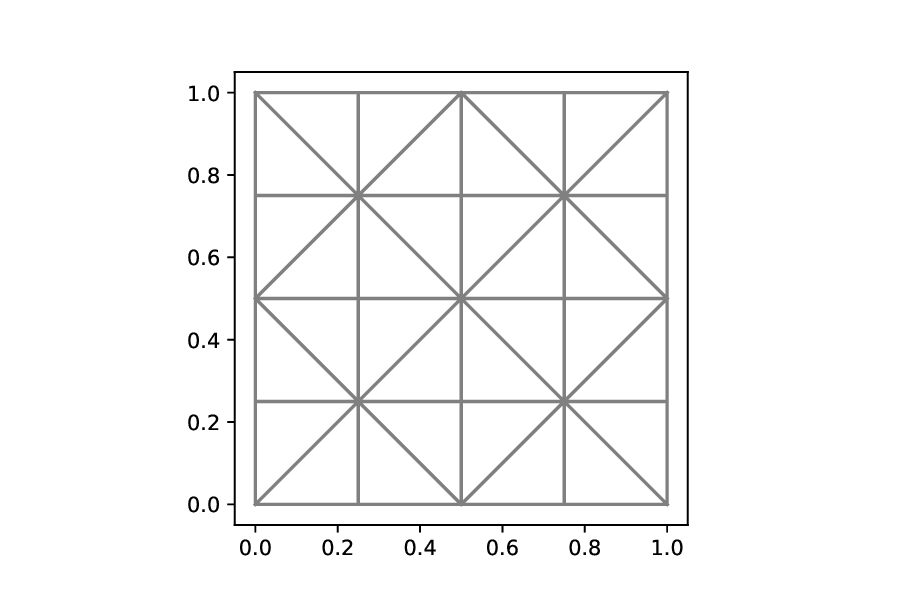}
\end{center}
\caption{\label{fig:uniform_mesh_square} The uniform mesh with $h=1/4$ for the unit square.%
}
\label{fi:squaremesh}
\end{figure}

\begin{table}
\caption{\label{ta:square_domain_clusters}The four leading clusters for the unit square.}
\begin{center}
\begin{tabular}{|c|l|}
\hline
\rule[-2mm]{0cm}{6mm}{}
Cluster & Eigenvalues   \\
\hline
1 & $\lambda_1 = 2\pi^2$  \\
\hline
2 & $\lambda_2 = \lambda_3 = 5\pi^2$   \\ 
\hline
3 & $\lambda_4 = 8\pi^2$   \\ 
\hline
4 & $\lambda_5 = \lambda_6 = 10\pi^2$  \\ 
\hline
\end{tabular}
\end{center}
\end{table}

For each cluster $K = 1, 2, 3, 4$, we compute bounds of  on $\delta_a(E_K,\hE_K)$ and $\delta_b(E_K,\hE_K)$ by applying Algorithm I 
(Theorem \ref{th:mainenergy} along with \eqref{eq:a_by_b}) and Algorithm II (Theorem \ref{thm:extension-davis-kahan-for-delta-b} along with  \eqref{eq:a_by_b}), respectively. 

The convergence behavior of computed bounds from Algorithm I 
on a sequence of uniformly refined meshes, when applied to the exact directed distances of the four leading clusters, is shown in Figures~\ref{fi:squareH1} and \ref{fi:squareL2}.
Particularly, the results confirm the expected optimal convergence rate $O(h)$ of the estimate \eqref{eq:Deltaest},
the convergence rate $O(h^2)$ of $\delta_b(E_K,\hE_K)$, and the sub-optimal rate $O(h)$ of the bound \eqref{eq:Deltaest_b}.

The corresponding indices of effectivity $\IndEff{\delta}$, for the the estimates computed by both algorithms are listed in Table~\ref{tab:square_domain_effectivity_index} and Table~\ref{table:index_residual_err_unit_square}.
The estimate by Algorithm II provides impressively sharp bounds for the error of approximate eigenfunctions.
Observed from the effectivity index, we can see that the estimate from Algorithm II on $\delta_b$ has an optimal convergence rate as $O(h^2)$.

Note that the estimate computed by Algorithm I only utilizes the Rayleigh quotient of the approximate eigenfunction, while Algorithm II takes advantage of the residual error estimation along with the flux constructed by solving variational equations.
The high precision of bounds on $\delta_a$ and $\tilde{\delta}$ based on
Theorem \ref{thm:extension-davis-kahan-for-delta-b} may be attributed to the relations \eqref{eq:a_by_b} and \eqref{eq:est-tilde-delta-by-delta-b}, which come from the fundamental relation in \eqref{eq:fundamental-relation-between-eigenfunction-and-eigenvalue}.
The a priori error estimates with the finite element projection can also help to obtain 
sharp estimates of $\delta_a$ and $\tilde{\delta}$ achieving effectivity indices almost one. This approach will be investigated in our subsequent paper. See also the discussion in preprint \cite[\S 6]{liu-Vejchodsky-arxiv-v1}.

\begin{figure}[tbh]
    \centering
    \includegraphics[width=\textwidth]{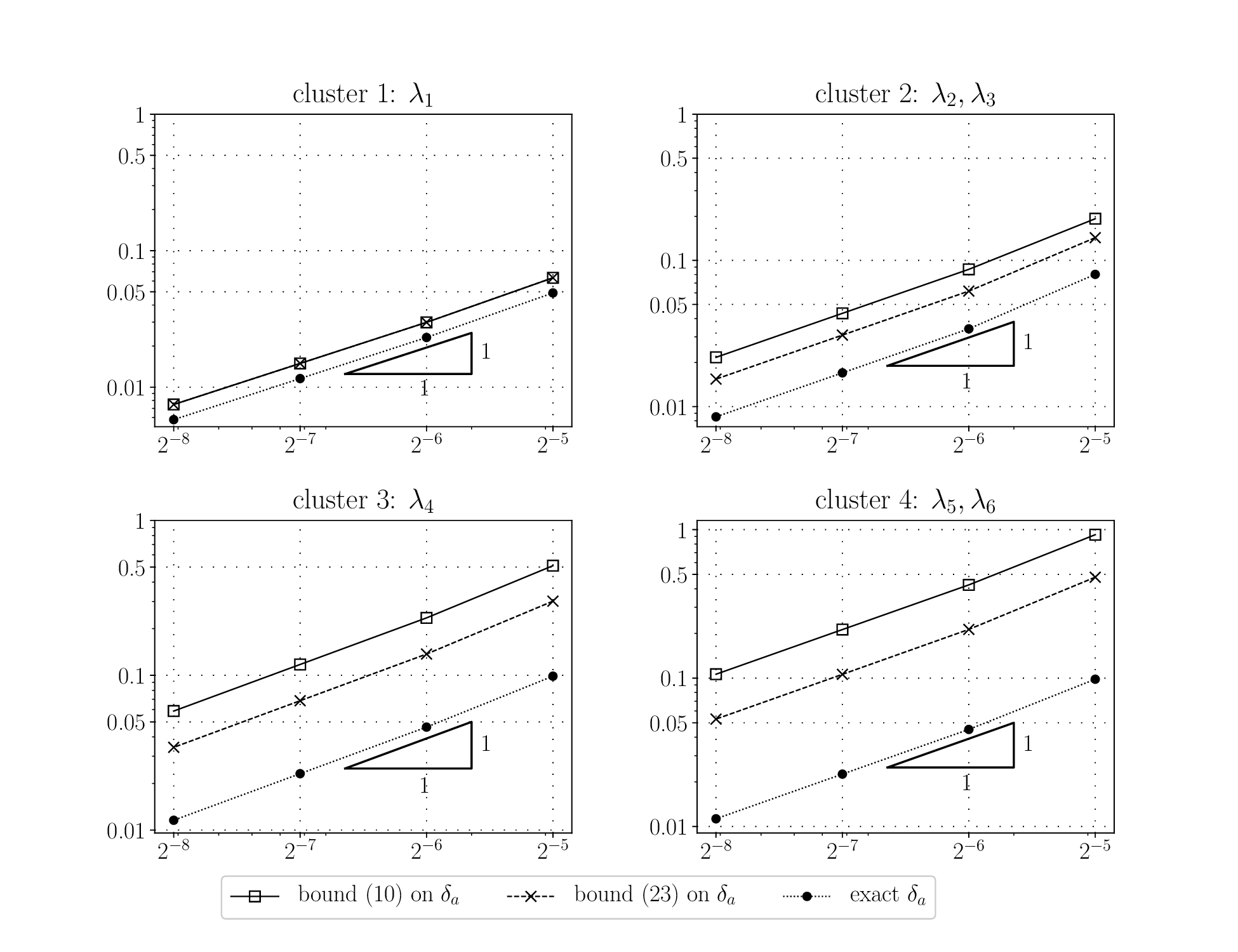}
    \caption{Bounds \eqref{eq:Deltaest} and \eqref{eq:a_by_b} on the energy distances of spaces of eigenfunctions $\delta_a(E_K,\hE_K)$ for the square domain and the four leading clusters of eigenvalues. 
    }
    \label{fi:squareH1}
\end{figure}
\begin{figure}[tbh]
    \centering
    \includegraphics[width=\textwidth]{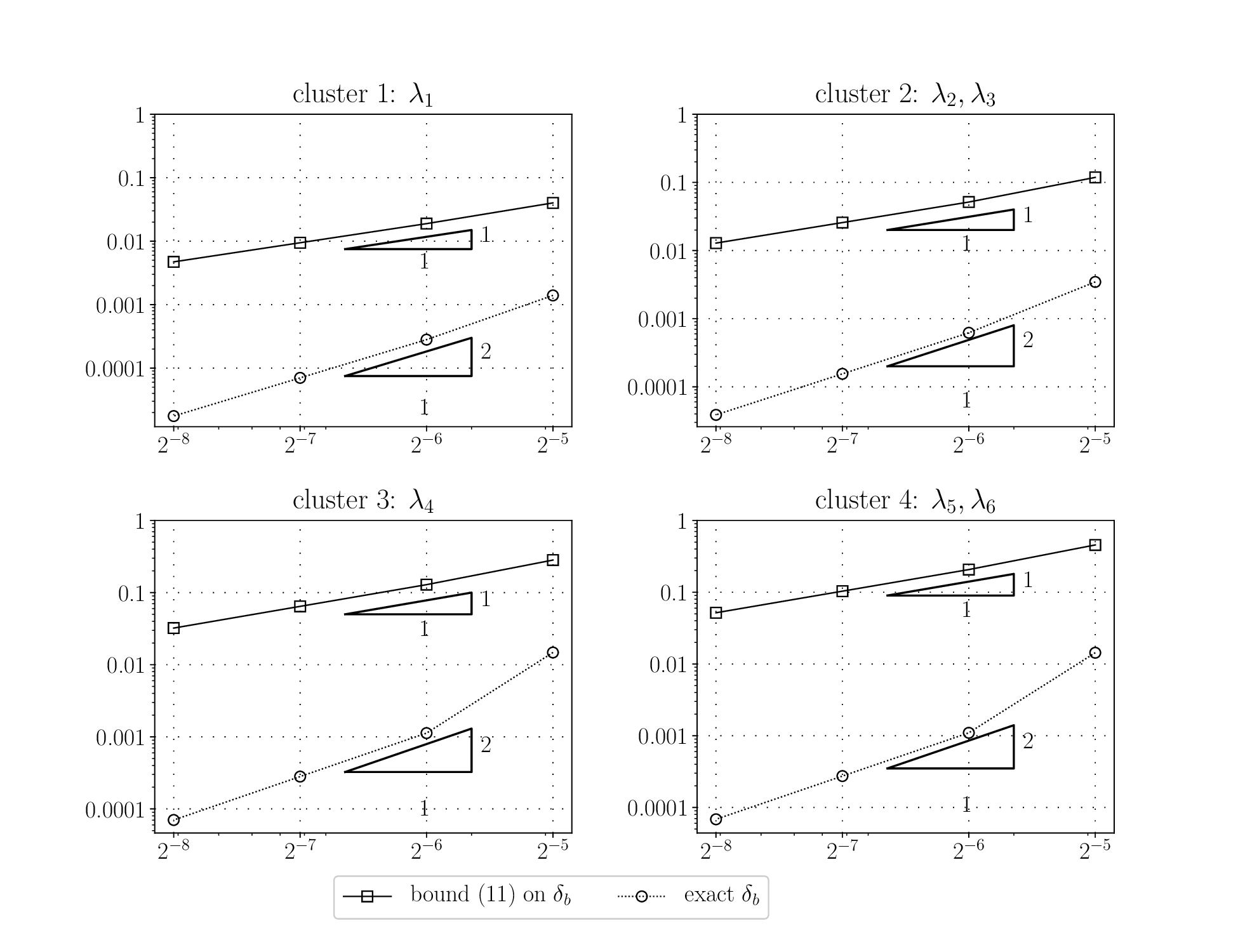}
    \caption{Bounds \eqref{eq:Deltaest_b} on the $L^2(\Omega)$ distances of spaces of eigenfunctions
    $\delta_b(E_K,\hE_K)$ for the square domain and four leading clusters of eigenvalues. 
    }
    \label{fi:squareL2}
\end{figure}

\begin{table}[ht]
    \caption{Effectivity indices of estimates \eqref{eq:Deltaest}, \eqref{eq:a_by_b} and \eqref{eq:Deltaest_b} (without the residual error estimation; the unit square)}
    \label{tab:square_domain_effectivity_index}
    \centering
    \begin{tabular}{|l|l|l|l|l|l|l|l|l|}
    \hline
\multirow{2}{*}{mesh size} &
\multicolumn{2}{l|}{\rule[-2mm]{0mm}{6mm}{}$\delta_a(E_1,\hE_1)$ } & \multicolumn{2}{l|}{$\delta_a(E_2,\hE_2)$ }& \multicolumn{2}{l|}{$\delta_a(E_3,\hE_3)$ } & \multicolumn{2}{l|}{$\delta_a(E_4,\hE_4)$ }  \\
     \cline{2-9} 
    & \rule[-1.5mm]{0mm}{5mm}{} \eqref{eq:Deltaest}  & \eqref{eq:a_by_b}  & \eqref{eq:Deltaest}  & \eqref{eq:a_by_b}  & \eqref{eq:Deltaest}  & \eqref{eq:a_by_b}  & \eqref{eq:Deltaest}  & \eqref{eq:a_by_b} \\
        \hline
    h=1/32 & 1.29 & 1.29 & 2.53 & 1.81 & 5.08 & 2.98 & 9.41 & 4.80 \\
        \hline
    h=1/64 &   1.29 & 1.29 & 2.55 & 1.81 & 5.08 & 2.96 & 9.42 & 4.72\\
        \hline
    h=1/128 &  1.29 & 1.29 & 2.55 & 1.81 & 5.08 & 2.96 & 9.42 & 4.70\\
        \hline
\end{tabular} \\
\vskip 0.2cm
    \begin{tabular}{|l|l|l|l|l|}
    \hline
\multirow{2}{*}{mesh size} & \rule[-2mm]{0mm}{6mm}{} $\delta_b(E_1,\hE_1)$  & $\delta_b(E_2,\hE_2)$  & $\delta_b(E_3,\hE_3)$  & $\delta_b(E_4,\hE_4)$ \\
     \cline{2-5} 
    & \rule[-1.5mm]{0mm}{5mm}{}  (10)  & (10) & (10)  & (10)   \\
        \hline
    h=1/32 &  33.57 & 41.28 & 56.90 & 93.58 \\
        \hline
    h=1/64 &  67.25 & 82.88 & 114.50 & 188.63 \\
        \hline
    h=1/128 & 134.56 & 165.91 & 229.34 & 378.00 \\
        \hline
\end{tabular}

\end{table}

\begin{table}[ht]
    \centering
    \caption{\label{table:index_residual_err_unit_square} Effectivity indices of estimation from Theorem \ref{thm:extension-davis-kahan-for-delta-b} (with the residual error estimation; the unit square)}
    \begin{tabular}{|l|l|l|l|l|l|l|l|l|}
    \hline
\multirow{2}{*}{mesh size} &
\multicolumn{2}{l|}{ \rule[-2mm]{0mm}{6mm}{}  cluster 1} & \multicolumn{2}{l|}{cluster 2}  & \multicolumn{2}{l|}{ cluster 3} &  \multicolumn{2}{l|}{cluster 4} \\
     \cline{2-9} 
    & \rule[-1.5mm]{0mm}{5mm}{}   $\IndEff{\delta_a}$ & $\IndEff{\delta_b}$  &  $\IndEff{\delta_a}$ & $\IndEff{\delta_b}$ &  $\IndEff{\delta_a}$ & $\IndEff{\delta_b}$  &  $\IndEff{\delta_a}$ & $\IndEff{\delta_b}$  \\
        \hline
    h=1/32  & 1.00045 & 1.27  &1.0065 & 3.12 & 1.078 & 8.189 & 1.12 & 10.385\\
        \hline
    h=1/64  & 1.00011 & 1.27 &1.0016 & 3.12 & 1.018 & 7.76& 1.026& 9.53 \\
        \hline
    h=1/128 & 1.00003 & 1.27 &1.0004 &  3.12 & 1.0043 & 7.66  & 1.0064 & 9.36\\
        \hline
\end{tabular}

\end{table}

Note that all bounds in this and in the following sections are computed in the floating-point arithmetic, and the influence of rounding errors is not taken into account. However, if needed, mathematically rigorous estimates could be obtained by employing the interval arithmetic \cite{moore2009introduction}.

\begin{remark}
Error estimates of eigenfunctions for the Dirichlet Laplacian in the square domain are well discussed in the existing literature.
For example, in \cite{CanDusMadStaVoh2017}, a computable estimate for the error in eigenfunctions is provided by solving an auxiliary problem and considering estimates of the residual, which are also used in Algorithm II. Note that paper \cite{CanDusMadStaVoh2017} measures the error of eigenfunctions with a different distance. 
Let us define the distance $\overline{\delta}(E, \hE )$ by
\begin{equation}
\label{eq:est-another-norm}
\overline{\delta}(E, \hE )
:=\max_{ \substack{\hat{u} \in \hE \\ \norm{\gamma \hat{u}}_b = 1} }  \min_{\substack{ {u} \in E \\ \norm{\gamma u}_b = 1}} \| u- \hat{u}\|_a.
\end{equation}
Note that $\overline{\delta}(E, \hE )=\overline{\delta}(\hE,E)$, 
only if $E$ and $\hE$ have the same dimension and the cluster width is zero. 
Paper \cite{CanDusMadStaVoh2017} utilizes this distance for the case of simple eigenvalues.
It is easy to see that $\tilde{\delta} \le \overline{\delta}$. By confirming the utilization of $\tilde{\delta}$ in Theorem \ref{lem:relation_of_tilde_delta_and_delta_b} and 
Theorem \ref{thm:extend-of-davis-kahan-for-delta-b},
the arguments therein also hold by replacing $\tilde{\delta}$ with $\overline{\delta}$. Thus, we can apply the bound for $\tilde{\delta}$ to directly  estimate $\overline{\delta}$.
%
That is, 
\begin{equation}
\label{eq:bar-delta-est}
\overline{\delta}^2(E, \hE)
\le 
\lambda_N + \hat{\lambda}_N - 2 \lambda_n \sqrt{1-\delta_b^2(E,\hE)}~.
\end{equation}
In  Table~\ref{tab:unit_square_vs_martin}, the estimate of $\overline{\delta}$ 
through \eqref{eq:bar-delta-est} and the estimate of $\delta_b$ from Theorem \ref{th:mainenergy} and Theorem \ref{thm:extension-davis-kahan-for-delta-b} is applied to 
the first cluster, and the effectivity indices of the obtained error estimates
are compared with the ones of \cite{CanDusMadStaVoh2017}. 
Note that in this comparison, we follow the mesh type used in \cite{CanDusMadStaVoh2017} rather than the type in Figure \ref{fig:uniform_mesh_square}
and the eigenvalue bounds for the clusters are different from the ones used in \cite{CanDusMadStaVoh2017}.
Further note the convergence rate $O(h)$ of the overestimation factor of \cite{CanDusMadStaVoh2017}. 
With the analogous argument as in Remark \ref{remark:super-convergence}, the locally reconstructed flux can also be used in \eqref{eq:estimate-for-espilon} to obtain sharp estimates with the quadratic convergence rate of the overestimating factor.
\end{remark}

\begin{table}[h]
    \centering
    \caption{Comparison of the effectivity index for different approaches applied to $\overline{\delta}(E_1, \hE_1)$ (the same type uniform mesh as \cite{CanDusMadStaVoh2017} for the unit square domain )}
    \begin{tabular}{|c|c|c|c|c|c|}
    \hline
    \rule[-2mm]{0mm}{6mm}{} $h$ & $\overline{\delta}$ & $\IndEff{\overline{\delta}}$ from Thm \ref{th:mainenergy} & $\IndEff{\overline{\delta}}$ from Thm \ref{thm:extension-davis-kahan-for-delta-b} &
    $\IndEff{\overline{\delta}}$ from 
    \eqref{eq:estimate-for-a-norm-by-davis-kahan-direct}
    &
    $\IndEff{\overline{\delta}}$ of \cite{CanDusMadStaVoh2017}  \\
        \hline
    \rule[-2mm]{0mm}{6mm}{}
    1/20 & $0.3494$  & $1.29$ & 1.00035 & 1.67 & 1.13 \\
        \hline
        \rule[-2mm]{0mm}{6mm}{}
        1/40 & 
        $0.1745$  & 1.29 & 1.00008 & 1.67 & 1.07 \\
        \hline
    \end{tabular}
    \label{tab:unit_square_vs_martin}
\end{table}

\subsection{The L-shaped domain}
We consider Laplace eigenvalue problem \eqref{eq:modpro} in the L-shaped domain $\Omega = (-1,1)^2\setminus(-1,0]^2$
to present the standard example with singularities of eigenfunctions
and also to demonstrate the versatility of the proposed method. We solve this problem by using both the classical linear conforming finite element space and the extended space with two more basis functions $$
\eta_k = (x^4-1)(y^4-1) r^{\frac{2k}{3}} \sin\left(\frac{2k}{3} (\theta + \pi/2) \right),\quad k=1,2.
$$

Since the exact eigenvalues are not known, bounds \eqref{eq:Deltaest},
\eqref{eq:Deltaest_b}, and \eqref{eq:a_by_b} 
are evaluated by using
two-sided bounds on eigenvalues, which were computed in \cite{liu2014high} 
and we list them in Table~\ref{tab:l_shaped_eig_lower_bound}. 
The first four eigenvalues are simple and form trivial clusters.

\begin{table}[h]
    \centering
    \caption{Lower bounds on the leading eigenvalues for the L-shaped domain.}
    \begin{tabular}{|c|c|c|c|c|}
    \hline
    \rule[-2mm]{0mm}{6mm}{}
    $\lambda_1$  & $\lambda_2$ &$\lambda_3$ & $\lambda_4$ & $\lambda_5$ \\
        \hline
    \rule[-2mm]{0mm}{6mm}{}
    $9.6397_{1}^{3}$ & $15.1972_{5}^{6}$ & $19.7392_0^1$ & $29.5214_7^9$ & $31.9126_2^4$ \\
        \hline
    \end{tabular}

    \label{tab:l_shaped_eig_lower_bound}
\end{table}

The initial finite element mesh is displayed in Figure~\ref{fig:l_shaped_domain}.
First, we apply Algorithm I to the four leading eigenvalue clusters.
Figure~\ref{fig:l_shape_singular_base_H1} shows the obtained bounds \eqref{eq:Deltaest} 
and \eqref{eq:a_by_b} 
on the energy distance $\delta_a$
and 
Figure~\ref{fig:l_shape_singular_base_L2} shows the bound \eqref{eq:Deltaest_b} on the $L^2(\Omega)$ distance $\delta_b$. 
The results confirm that both the distance $\delta_a$ and its bounds converge with a slower rate for the standard finite element space and have the optimal speed of convergence $O(h)$ for the extended finite element space.

The indices of effectivity for the estimates computed by the two algorithms are listed in Table~\ref{table:l-shape-eff-index} and the comparison to the result of \cite{CanDusMadStaVoh2017} is shown in  Table~\ref{table:l-shape-comparison-with-martin}. Note the exact eigenfunctions are unknown and the approximate ones over quite refined mesh are used as the exact eigenfunctions in the effectivity index computation. The stable values of the 
effectivity indices indicate that the bound computed by 
Algorithm II has the optimal rate of convergence. 

\begin{figure}
    \centering
    \includegraphics[scale=0.3]{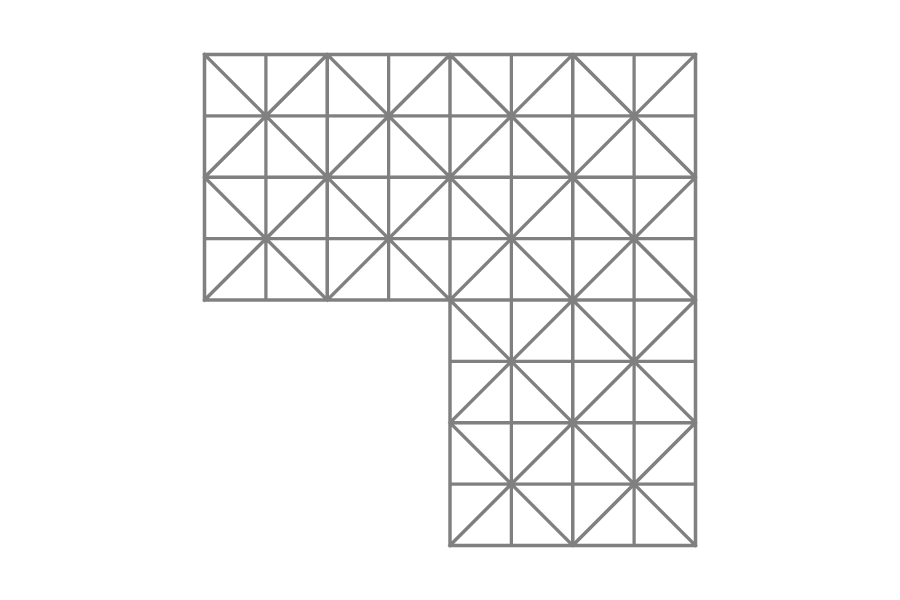}
    \caption{L-shaped domain and the initial mesh
    }
    \label{fig:l_shaped_domain}
\end{figure}

\begin{figure}[tbh]
    \centering
\includegraphics[width=\textwidth]{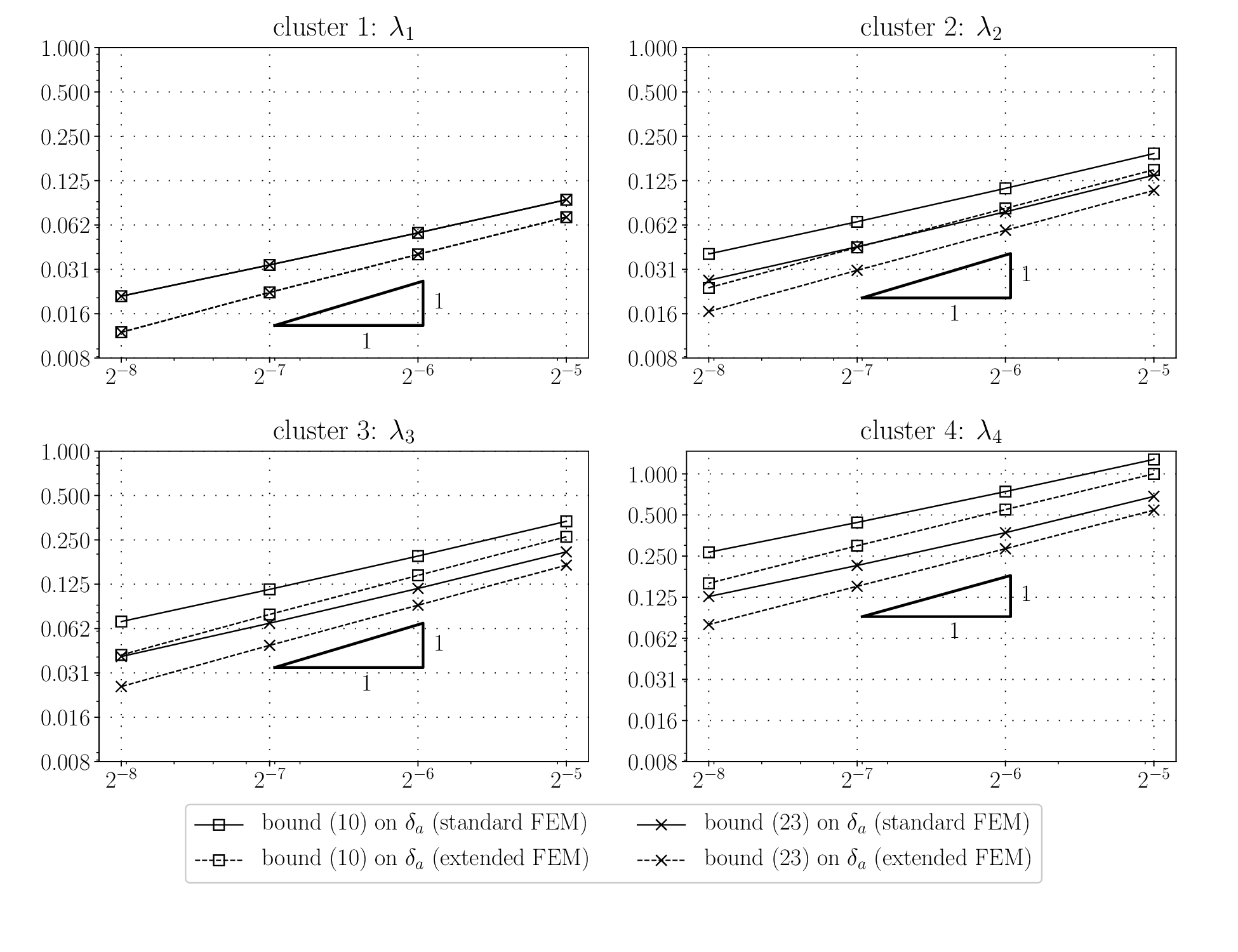}
    \caption{%
    Bounds \eqref{eq:Deltaest} and \eqref{eq:a_by_b} on the energy distance $\delta_a$  of spaces of eigenfunctions for the L-shaped domain.
    Results for the standard finite element space are plotted with solid lines while results for the extended space with dashed lines.
\label{fig:l_shape_singular_base_H1}
    }
\end{figure}

\begin{figure}[tbh]
    \centering
\includegraphics[width=\textwidth]{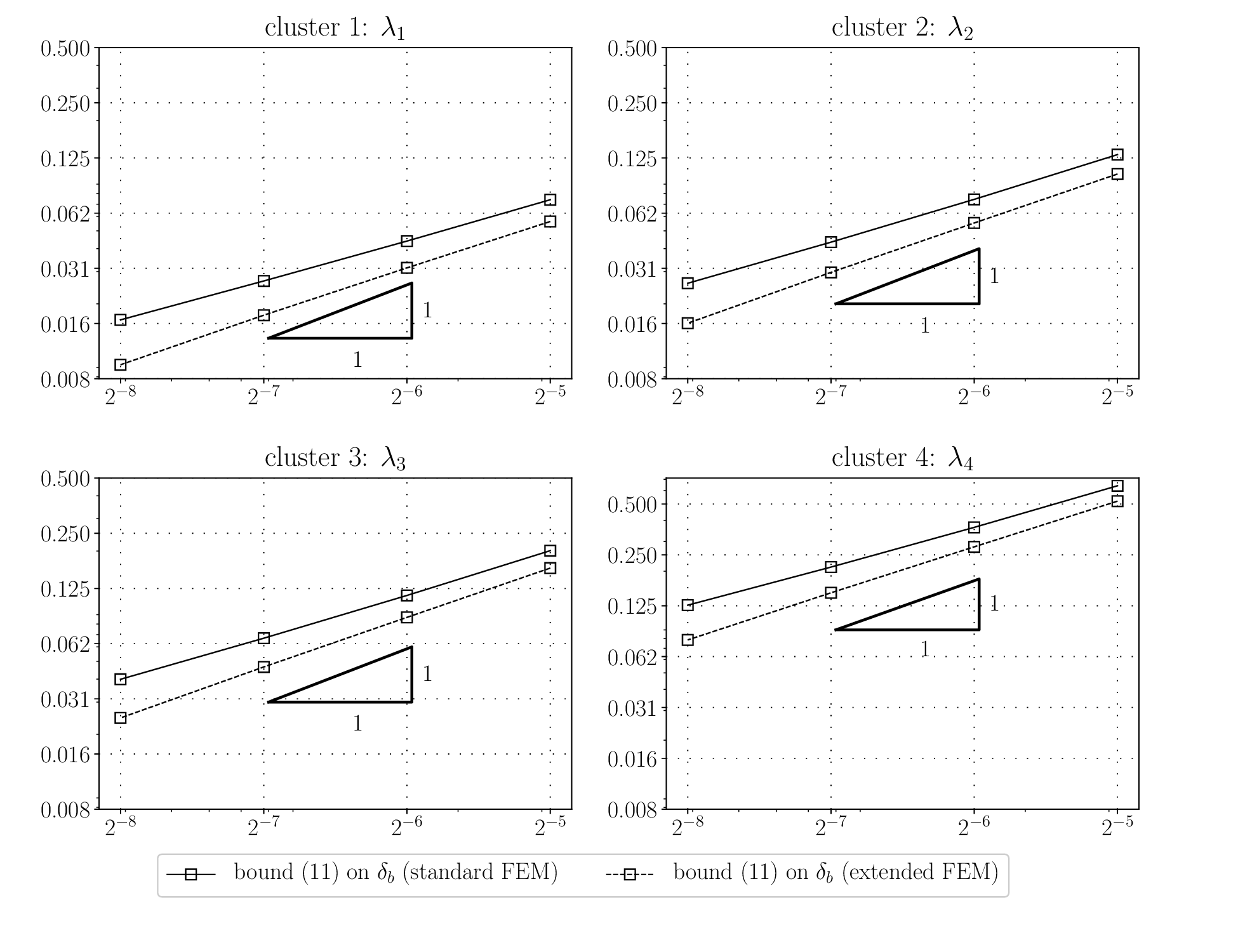}%
\vskip -0.5cm
    \caption{%
    Bounds \eqref{eq:Deltaest_b} on the $L^2(\Omega)$ distance $\delta_b$ of spaces of eigenfunctions for the L-shaped domain.
    Results for the standard finite element space are plotted with solid lines while results for the extended space with dashed lines.
    \label{fig:l_shape_singular_base_L2}
    }
\end{figure}

\begin{table}[h]
    \centering
    \caption{\label{table:l-shape-eff-index} Effectivity index of estimation on $\delta_a$ and $\delta_b$ (standard FEM; L-shaped domain)}
    \begin{tabular}{|l|l|l|l|l|l|l|l|l|}
    \multicolumn{9}{c}{ \rule[-3mm]{0mm}{8mm}{}(a) Estimate of $\delta_a$ and $\delta_b$ by  \eqref{eq:a_by_b} and \eqref{eq:Deltaest_b}   (no residual error estimation)}\\
        \hline
\multirow{2}{*}{mesh size} &
\multicolumn{2}{l|}{ \rule[-2mm]{0mm}{6mm}{}  cluster 1} & \multicolumn{2}{l|}{cluster 2}  & \multicolumn{2}{l|}{ cluster 3} &  \multicolumn{2}{l|}{cluster 4} \\
     \cline{2-9} 
    & \rule[-1.5mm]{0mm}{5mm}{}   $\IndEff{\delta_a}$ & $\IndEff{\delta_b}$  &  $\IndEff{\delta_a}$ & $\IndEff{\delta_b}$  &  $\IndEff{\delta_a}$ & $\IndEff{\delta_b}$   &  $\IndEff{\delta_a}$ & $\IndEff{\delta_b}$   \\
        \hline
    h=1/32  & 1.67 & 24.99  & 3.57 & 159.94 & 4.48 & 178.45 & 12.74 & 394.83\\
        \hline
    h=1/64  & 1.69  & 38.71 & 4.03 & 365.49 & 5.056 & 407.52& 13.83 & 898.22 \\
        \hline
\end{tabular}
\begin{tabular}{|l|l|l|l|l|l|l|l|l|}
    \multicolumn{9}{c}{ \rule[-3mm]{0mm}{8mm}{} (b) Estimate of $\delta_a$ and $\delta_b$ by Theorem \ref{thm:extension-davis-kahan-for-delta-b} (residual error estimation used)}\\
    \hline
\multirow{2}{*}{mesh size} &
\multicolumn{2}{l|}{ \rule[-2mm]{0mm}{6mm}{}  cluster 1} & \multicolumn{2}{l|}{cluster 2}  & \multicolumn{2}{l|}{ cluster 3} &  \multicolumn{2}{l|}{cluster 4} \\
     \cline{2-9} 
    & \rule[-1.5mm]{0mm}{5mm}{}    $\IndEff{\delta_a}$ & $\IndEff{\delta_b}$  &  $\IndEff{\delta_a}$ & $\IndEff{\delta_b}$  &  $\IndEff{\delta_a}$ & $\IndEff{\delta_b}$   &  $\IndEff{\delta_a}$ & $\IndEff{\delta_b}$    \\
        \hline
    h=1/32  & 1.012 & 1.97  & 1.0087 & 5.99 & 1.021 & 8.45 & 1.34 & 29.36\\
        \hline
    h=1/64  & 1.024 & 1.82 & 1.0021 & 5.98 & 1.0052 & 8.33& 1.060& 23.26 \\
        \hline
\end{tabular}
\end{table}

\begin{table}[h]
    \centering
    \caption{\label{table:l-shape-comparison-with-martin} Comparison of the effectivity index for different approaches (L-shaped domain)}
    \begin{tabular}{|c|c|c|c|}
    \hline
    \rule[-2mm]{0mm}{6mm}{Approaches}  & Mesh type & \#(DOF) &
    $\IndEff{\overline{\delta}(E_1, \hE_1)}$  \\
        \hline
    \rule[-2mm]{0mm}{6mm}{} Theorem \ref{th:mainenergy}  & uniform mesh ($h=1/64$) & 12545 & 1.69 \\
        \hline
    \rule[-2mm]{0mm}{6mm}{} Theorem \ref{thm:extension-davis-kahan-for-delta-b}  & uniform mesh ($h=1/64$) & 12545 & 1.024 \\
    \hline
    \rule[-2mm]{0mm}{6mm}{} Estimate of \cite{CanDusMadStaVoh2017}  & unstructured mesh & 24925 & 2.51  \\
        \hline
    \end{tabular}
    \label{tab:L_shape_vs_martin}
\end{table}
\subsection{The dumbbell shaped domain}

\begin{figure}[tbh]
    \centering
    \includegraphics[height=3cm]{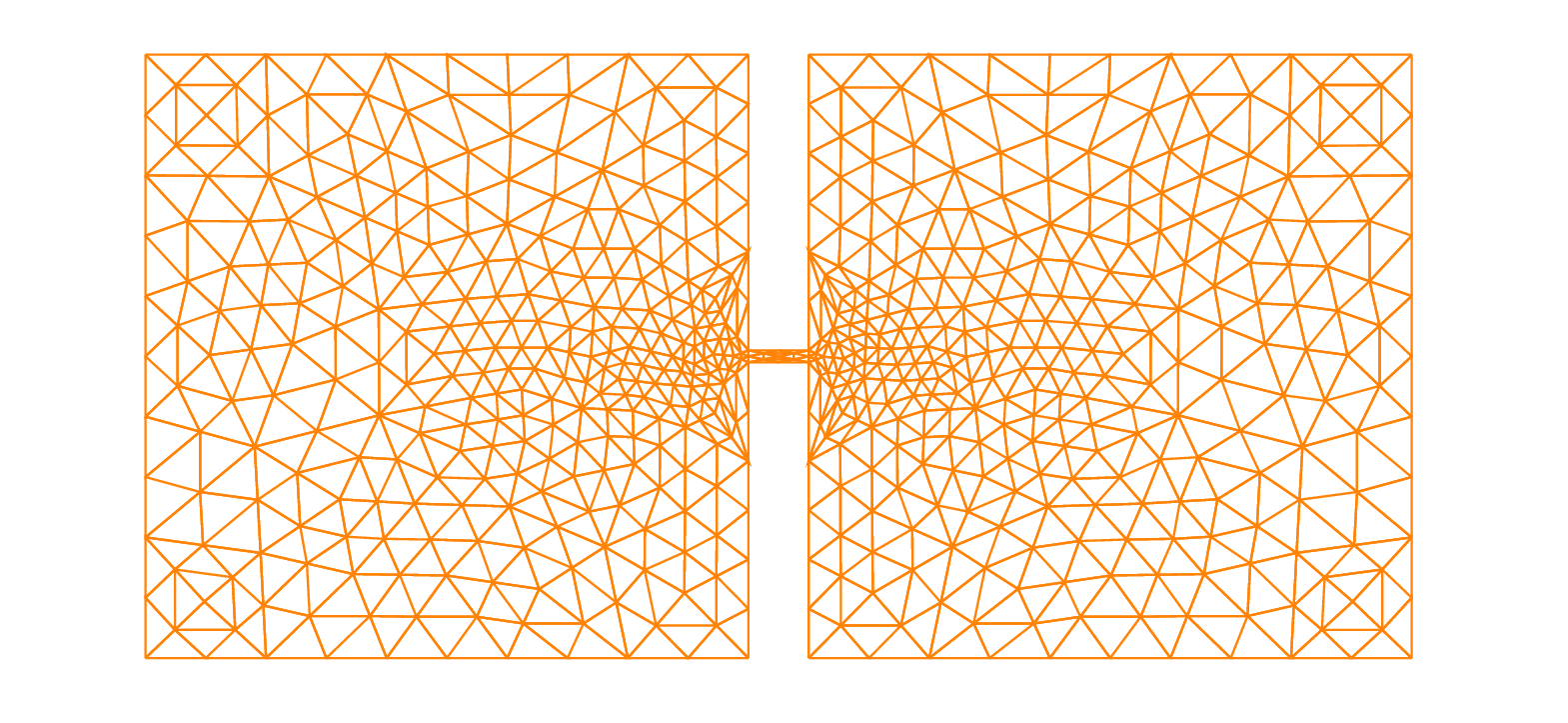}
    \caption{Dumbbell-shaped domain and the initial mesh}
    \label{fi:dumbbell_domain}
\end{figure}

\begin{table}[ht]
\begin{center}
\caption{\label{ta:dumbbell_domain_clusters} Lower and upper bounds of eigenvalues for the dumbbell shaped domain. Two times refined initial mesh and third order finite element spaces were used.}
\begin{tabular}{|c|l|}
\hline
    \rule[-2mm]{0mm}{6mm}{}
cluster& lower and upper eigenvalue bounds \\
\hline
    \rule[-2mm]{0mm}{6mm}{}
1 & $\lambda_1=19.736^{729}_{634},~ \lambda_2=19.736^{729}_{635}$  \\
\hline
    \rule[-2mm]{0mm}{6mm}{}
2 & $\lambda_3=49.33^{809}_{761},~\lambda_4=49.33^{809}_{761},~ \lambda_5=49.348020_{5}^{8},~\lambda_6=49.348020_{5}^{8}$ \\ 
\hline
    \rule[-2mm]{0mm}{6mm}{}
3 & $\lambda_7=78.9568_{290}^{301},~ \lambda_8=78.9568_{290}^{301}$ \\ 
\hline
    \rule[-2mm]{0mm}{6mm}{}
4 & $\lambda_9=98.6_{69041}^{71154}, ~ \lambda_{10}=98.6_{69041}^{71154}, ~ \lambda_{11}=98.69604_{39}^{41}, ~\lambda_{12}=98.69604_{39}^{41}$ \\ 
\hline
\end{tabular}
\end{center}
\end{table}

Finally, we consider Laplace eigenvalue problem \eqref{eq:modpro}
in a dumbbell shaped domain consisting of two unit squares connected by a bar 
of width $0.02$ and length $0.1$, see Figure~\ref{fi:dumbbell_domain}, 
where also the initial mesh is depicted.
This example is interesting due to singularities of eigenfunctions in re-entrant corners and especially due to tight clusters of eigenvalues.

The exact eigenpairs are not known, but the eigenvalues are expected 
to be close to eigenvalues for the union of two squares, i.e., two eigenvalues close to 
$2\pi^2 \approx 19.739$, four eigenvalues close to $5\pi^2 \approx 49.348$, etc.
In order to compute high precision two-sided bounds for these eigenvalues, we
combine the Crouzeix--Raviart nonconforming finite elements and 
the Lehmann--Goerisch method as proposed in \cite{Liu2015}.
The resulting two-sided bounds obtained on a fine mesh and finite element spaces of the third order
are presented in Table~\ref{ta:dumbbell_domain_clusters}. 

Table~\ref{ta:dumbbell_domain_clusters} also shows the chosen division of the first twelve eigenvalues into four clusters. 
Note that eigenvalues $\lambda_3$ and $\lambda_4$ are strictly separated from $\lambda_5$ and $\lambda_6$. Therefore, they could be considered as two separate clusters, but then the spectral gap between them would be small and the factor $\rho - \lambda_n$ in \eqref{eq:Deltaest} and \eqref{eq:Deltaest_b} would yield large overestimation. For this reason, all four eigenvalues $\lambda_3, \dots, \lambda_6$ are considered in one cluster.

We apply compute the bounds on $\delta_a(E_K,\hE_K)$ and $\delta_b(E_K,\hE_K)$ for the four clusters $K=1,2,3,4$ as we did for the unit square domain. 
To investigate the rate of convergence of estimate computed by Algorithm I, 
the bounds \eqref{eq:Deltaest}, \eqref{eq:Deltaest_b}, and \eqref{eq:a_by_b}
are displayed in Figure~\ref{fi:dumbbellH1}.
The first and the third cluster are very tight and we observe the first order convergence. However, the convergence curves for the second and the fourth cluster bend due to the non-negligible width of these clusters.

The comparison of effectivity indices of the estimates from Algorithm I and II are shown in Table \ref{table:dumbbell-shape-eff-index}.
Here, the solution solved by the conforming Lagrange FEM with degree 4 is regarded as the exact solution in calculating the effectivity index. The mesh size $h$ here denotes the maximum edge length of the mesh. 
The presented results show that Algorithm II considerably improves the accuracy of the bounds, but has higher computational cost. %
{The computational cost can be decreased by computing only approximate minimum in \eqref{eq:estimate-for-espilon-simple} using local flux reconstructions as in \cite{CanDusMadStaVoh2017,CanDusMadStaVoh2019}. See also study \cite{Vejchodsky2018b} showing that the local flux reconstructions are considerably faster to compute than the global ones, while the two approaches provide almost the same accuracy.
}
\begin{figure}[htb]
    \centering
    \includegraphics[width=1.1\textwidth]{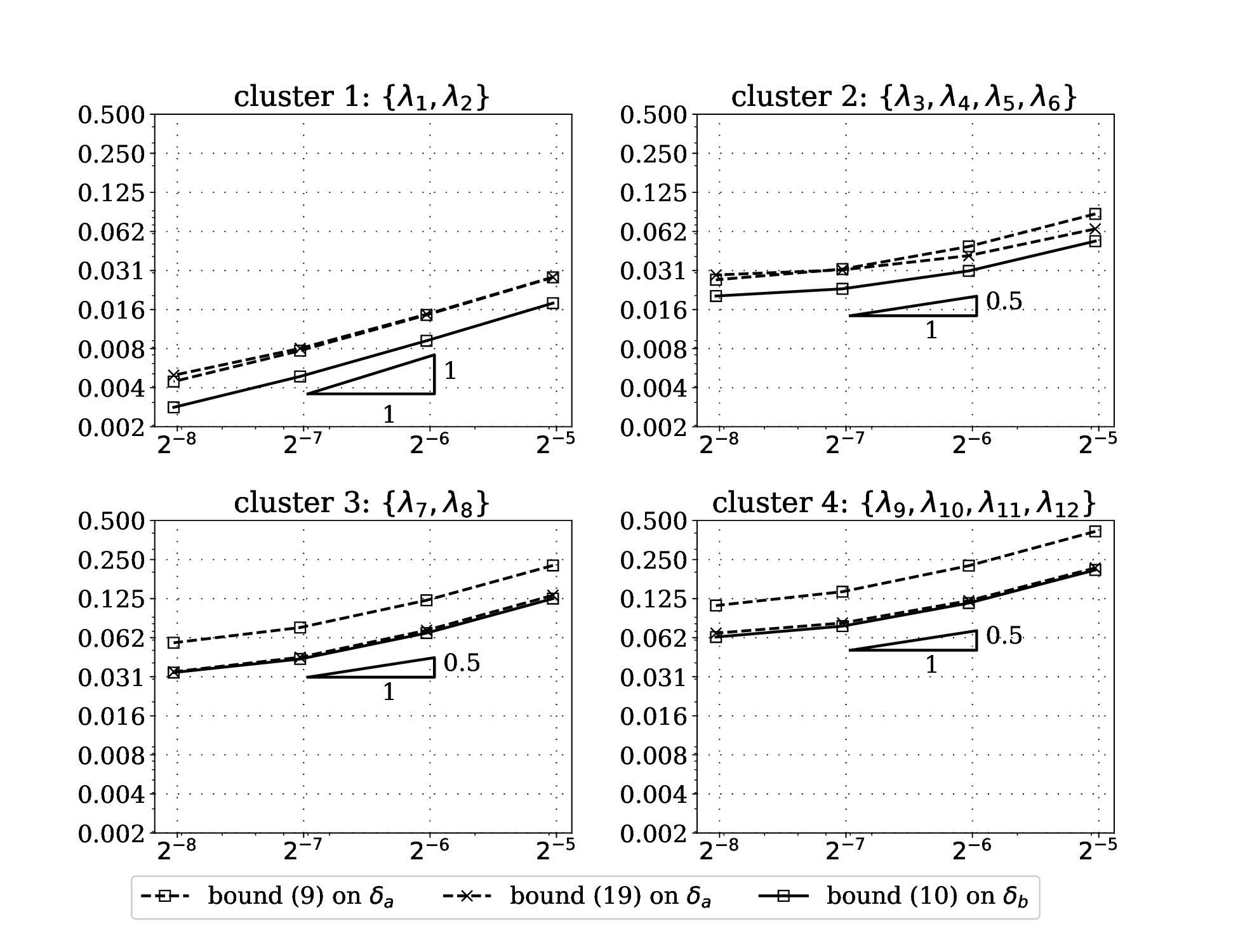}
    \caption{Bounds \eqref{eq:Deltaest}, \eqref{eq:Deltaest_b}, and \eqref{eq:a_by_b} 
    on the energy distance $\delta_a$ and $L^2(\Omega)$ distance $\delta_b$ of 
    spaces of eigenfunctions for the dumbbell shaped domain.
    }
    \label{fi:dumbbellH1}
\end{figure}

\begin{table}[h]
    \centering
    \caption{\label{table:dumbbell-shape-eff-index} Effectivity index of estimation on $\delta_a$ and $\delta_b$ (Dumbbell-shaped domain)}
    \begin{tabular}{|l|l|l|l|l|l|l|l|l|}
    \multicolumn{9}{c}{ \rule[-3mm]{0mm}{8mm}{}(a) Estimate of $\delta_a$ and $\delta_b$ by  \eqref{eq:a_by_b} and \eqref{eq:Deltaest_b}  (no residual error estimation)}\\
        \hline
\multirow{2}{*}{mesh size} &
\multicolumn{2}{l|}{ \rule[-2mm]{0mm}{6mm}{}  cluster 1} & \multicolumn{2}{l|}{cluster 2}  & \multicolumn{2}{l|}{ cluster 3} &  \multicolumn{2}{l|}{cluster 4} \\
     \cline{2-9} 
    & \rule[-1.5mm]{0mm}{5mm}{}   $\IndEff{\delta_a}$ & $\IndEff{\delta_b}$  &  $\IndEff{\delta_a}$ & $\IndEff{\delta_b}$  &  $\IndEff{\delta_a}$ & $\IndEff{\delta_b}$   &  $\IndEff{\delta_a}$ & $\IndEff{\delta_b}$   \\
        \hline
    h=0.068  & 1.29 & 29.94  & 2.41 & 29.64 & 5.10 & 40.27 & 8.32 & 42.96\\
        \hline
\end{tabular}
\begin{tabular}{|l|l|l|l|l|l|l|l|l|}
    \multicolumn{9}{c}{ \rule[-3mm]{0mm}{8mm}{} (b) Estimate of $\delta_a$ and $\delta_b$ by Theorem \ref{thm:extension-davis-kahan-for-delta-b} (residual error estimation used)}\\
    \hline
\multirow{2}{*}{mesh size} &
\multicolumn{2}{l|}{ \rule[-2mm]{0mm}{6mm}{}  cluster 1} & \multicolumn{2}{l|}{cluster 2}  & \multicolumn{2}{l|}{ cluster 3} &  \multicolumn{2}{l|}{cluster 4} \\
     \cline{2-9} 
    & \rule[-1.5mm]{0mm}{5mm}{}    $\IndEff{\delta_a}$ & $\IndEff{\delta_b}$  &  $\IndEff{\delta_a}$ & $\IndEff{\delta_b}$  &  $\IndEff{\delta_a}$ & $\IndEff{\delta_b}$   &  $\IndEff{\delta_a}$ & $\IndEff{\delta_b}$    \\
        \hline
    h=0.068  & 1.00033 & 1.044  & 1.028 & 2.38 & 1.064 & 5.26 & 1.17 & 5.84\\
    \hline
\end{tabular}
\end{table}

\section{Conclusions}
\label{se:conclusions}

For the abstractly formulated eigenvalue problem of compact operators over Hilbert spaces, we proposed two algorithms to derive guaranteed upper bounds 
on the directed distance between sub-spaces of exact and approximate eigenfunctions. The two presented algorithms have their own advantages for eigenvalue problems of different settings. Algorithm I only utilizes the Rayleigh quotients of the approximate eigenvectors, while Algorithm II provides sharper bounds by further considering the residual error estimation of discretized operators. 

The bounds of both algorithms are independent of the discretization method and apply to arbitrary conforming approximations of eigenfunctions.  
Bounds on the total error of these approximations are
easily computed by using solely the two-sided bounds on exact eigenvalues 
and the approximate eigenfunctions themselves.
The derived bounds can be straightforwardly applied to, for example, the (generalized) matrix, Laplace, Steklov, and many other eigenvalue problems. 

Numerical examples for the Laplace eigenvalue problem discretized by the finite element method show that the bound from Algorithm I in the $L^2$ norm converges with a sub-optimal rate, while the one from Algorithm II has the optimal rate.
In a subsequent work we would like to employ the Aubin--Nitsche technique and the explicit \emph{a priori} error estimate for the energy projection \cite{LiuOis2013} in order to derive bounds in the $L^2$ norm with the optimal rate of convergence, without solving the dual problem as required in Algorithm II.

\medskip

\section*{Conflict of interest}
The authors declare that they have no conflict of interest.

\medskip

\noindent
{\bf Acknowledgement } 
The authors greatly appreciate the valuable referees' comments. Thanks to them the original manuscript considerably improved.
The first author also show thanks to 
Dr. Yuji Nakatsukasa from Oxford University for his introduction of Davis--Kahan's method.


\providecommand{\bysame}{\leavevmode\hbox to3em{\hrulefill}\thinspace}
\providecommand{\MR}{\relax\ifhmode\unskip\space\fi MR }
\providecommand{\MRhref}[2]{%
  \href{http://www.ams.org/mathscinet-getitem?mr=#1}{#2}
}
\providecommand{\href}[2]{#2}

\bibliographystyle{spmpsci}

\begin{thebibliography}{10}

  \bibitem{ArmDur2004}
  Mar{\'{\i}}a~G. Armentano and Ricardo~G. Dur{\'a}n, \emph{Asymptotic lower
    bounds for eigenvalues by nonconforming finite element methods}, Electron.
    Trans. Numer. Anal. \textbf{17} (2004), 93--101 (electronic). \MR{2040799}
  
  \bibitem{BabOsb:1991}
  Ivo Babu{\v{s}}ka and John~E. Osborn, \emph{Eigenvalue problems}, Handbook of
    numerical analysis, {V}ol.\ {II}, North-Holland, Amsterdam, 1991,
    pp.~641--787. \MR{1115240}
  
  \bibitem{Behnke1991}
  Henning Behnke, \emph{The calculation of guaranteed bounds for eigenvalues using
    complementary variational principles}, Computing \textbf{47} (1991), no.~1,
    11--27. \MR{1137071}
  
  \bibitem{BirBooSwaWen:1966}
  Garrett Birkhoff, C.~de~Boor, B.~Swartz, and B.~Wendroff, \emph{Rayleigh-{R}itz
    approximation by piecewise cubic polynomials}, SIAM J. Numer. Anal.
    \textbf{3} (1966), 188--203. \MR{0203926}
  
  \bibitem{Boffi:2010}
  Daniele Boffi, \emph{Finite element approximation of eigenvalue problems}, Acta
    Numer. \textbf{19} (2010), 1--120. \MR{2652780 (2011e:65256)}
  
  \bibitem{CanDusMadStaVoh2017}
  Eric Canc{\`e}s, Genevi{\`e}ve Dusson, Yvon Maday, Benjamin Stamm, and Martin
    Vohral{\'\i}k, \emph{Guaranteed and robust a posteriori bounds for {L}aplace
    eigenvalues and eigenvectors: conforming approximations}, SIAM J. Numer.
    Anal. \textbf{55} (2017), no.~5, 2228--2254.
  
  \bibitem{CanDusMadStaVoh2018}
  \bysame, \emph{Guaranteed and robust a posteriori bounds for {L}aplace
    eigenvalues and eigenvectors: a unified framework}, Numer. Math. \textbf{140}
    (2018), no.~4, 1033--1079.
  
  \bibitem{CanDusMadStaVoh2019}
  \bysame, \emph{Guaranteed a posteriori bounds for eigenvalues and eigenvectors:
    multiplicities and clusters}, 
    Math. Comp. \textbf{89} (2020), no.~326, 2563--2611.

  \bibitem{CarGal2014}
  Carsten Carstensen and Dietmar Gallistl, \emph{Guaranteed lower eigenvalue
    bounds for the biharmonic equation}, Numer. Math. \textbf{126} (2014), no.~1,
    33--51. \MR{3149071}
  
  \bibitem{CarGed2014}
  Carsten Carstensen and Joscha Gedicke, \emph{Guaranteed lower bounds for
    eigenvalues}, Math. Comp. \textbf{83} (2014), no.~290, 2605--2629.
    \MR{3246802}
  
  \bibitem{Chatelin1983}
  Fran{\c{c}}oise Chatelin, \emph{Spectral approximation of linear operators},
    Academic Press, Inc., New York, 1983. \MR{716134}

  \bibitem{DarDurPad2012}
  Enzo Alberto Dari, Ricardo G. Dur\'{a}n, and Claudio Padra, \emph{A posteriori error estimates
    for non-conforming approximation of eigenvalue problems}, Appl. Numer. Math.
    \textbf{62} (2012), no.~5, 580--591. \MR{2899264}
  
  \bibitem{davis1970rotation}
  Chandler Davis and William~Morton Kahan, \emph{{The rotation of eigenvectors by
    a perturbation. III}}, SIAM J. Numer. Anal. \textbf{7} (1970), no.~1, 1--46.
  
  \bibitem{debnath2005introduction}
  Lokenath Debnath and Piotr Mikusinski, \emph{Introduction to {H}ilbert spaces
    with applications}, Academic press, 2005.
  
  \bibitem{DurGasPad1999}
  Ricardo~G. Dur\'{a}n, Lucia Gastaldi, and Claudio Padra, \emph{A posteriori
    error estimators for mixed approximations of eigenvalue problems}, Math.
    Models Methods Appl. Sci. \textbf{9} (1999), no.~8, 1165--1178. \MR{1722056}
  
  \bibitem{giani2018posteriori}
  Stefano Giani, Luka Grubi{\v{s}}i{\'c}, Harri Hakula, and Jeffrey~S Ovall,
    \emph{An a posteriori estimator of eigenvalue/eigenvector error for
    penalty-type discontinuous Galerkin methods}, Applied Mathematics and
    Computation \textbf{319} (2018), 562--574.
  
  \bibitem{GiaHal2012}
  Stefano Giani and Edward J.~C. Hall, \emph{An {\it a posteriori} error
    estimator for {$hp$}-adaptive discontinuous {G}alerkin methods for elliptic
    eigenvalue problems}, Math. Models Methods Appl. Sci. \textbf{22} (2012),
    no.~10, 1250030, 35~p. \MR{2974168}
  
  \bibitem{GoeHau1985}
  Friedrich Goerisch and Heinz Haunhorst, \emph{Eigenwertschranken f\"ur
    {E}igenwertaufgaben mit partiellen {D}ifferentialgleichungen}, Z. Angew.
    Math. Mech. \textbf{65} (1985), no.~3, 129--135. \MR{789949}
  
  \bibitem{HongXieYueZhang2018}
  Qichen Hong, Hehu Xie, Meiling Yue, and Ning Zhang, \emph{Fully computable
    error bounds for eigenvalue problem}, Int. J. Numer. Anal. Model. \textbf{15}
    (2018), no.~1-2, 260--270. \MR{3722957}
  
  \bibitem{HuHuaLin2014}
  Jun Hu, Yunqing Huang, and Qun Lin, \emph{Lower bounds for eigenvalues of
    elliptic operators: by nonconforming finite element methods}, J. Sci. Comput.
    \textbf{61} (2014), no.~1, 196--221. \MR{3254372}
  
  \bibitem{JiaCheXie2013}
  Shanghui Jia, Hongtao Chen, and Hehu Xie, \emph{A posteriori error estimator
    for eigenvalue problems by mixed finite element method}, Sci. China Math.
    \textbf{56} (2013), no.~5, 887--900. \MR{3047040}
  
  \bibitem{Kato1949}
  Tosio Kato, \emph{On the upper and lower bounds of eigenvalues}, J. Phys. Soc.
    Japan \textbf{4} (1949), 334--339. \MR{0038738}
  
  \bibitem{Lehmann1949}
  Nikolaus Joachim Lehmann, \emph{Beitr\"age zur numerischen {L}\"osung linearer
    {E}igenwertprobleme. {I}}, Z. Angew. Math. Mech. \textbf{29} (1949),
    341--356. \MR{0034511}
  
  \bibitem{Lehmann1950}
  \bysame, \emph{Beitr\"age zur numerischen {L}\"osung linearer
    {E}igenwertprobleme. {II}}, Z. Angew. Math. Mech. \textbf{30} (1950), 1--16.
    \MR{0034512}
  
  \bibitem{LiaoYuLiu2019}
  Shih-Kang Liao, Yu-Chen Shu, and Xuefeng Liu, \emph{Optimal estimation for the
    Fujino--Morley interpolation error constants}, Jpn. J. Ind. Appl. Math.
    \textbf{36} (2019), no.~2, 521--542.
  
  \bibitem{Liu2015}
  Xuefeng Liu, \emph{A framework of verified eigenvalue bounds for self-adjoint
    differential operators}, Appl. Math. Comput. \textbf{267} (2015), 341--355.
    \MR{3399052}
    
  \bibitem{Liu2020}
  Xuefeng Liu, \emph{Explicit eigenvalue bounds of differential operators defined by symmetric positive semi-definite bilinear forms}, 
  J. of Comp. \& Appl.  Math.,
  \textbf{371} (2020), 112666.
  
  \bibitem{LiuOis2013}
  Xuefeng Liu and Shin'ichi Oishi, \emph{Verified eigenvalue evaluation for the
    {L}aplacian over polygonal domains of arbitrary shape}, SIAM J. Numer. Anal.
    \textbf{51} (2013), no.~3, 1634--1654. \MR{3061473}
  
  \bibitem{liu2014high}
  Xuefeng Liu, Tomoaki Okayama, and Shin'ichi Oishi, \emph{{High-Precision
    Eigenvalue Bound for the Laplacian with Singularities}}, Computer
    Mathematics, Springer, 2014, pp.~311--323.
  
  \bibitem{liu-you:2018}
  Xuefeng Liu and Chun'guang You, \emph{{Explicit bound for quadratic Lagrange
    interpolation constant on triangular finite elements}}, Appl. Math. Comput.
    \textbf{319} (2018), 693--701.

  \bibitem{liu-Vejchodsky-arxiv-v1}
  Xuefeng Liu and Tom{\'a}{\v{s}} Vejchodsk{\'y},
  \emph{Rigorous and fully computable a posteriori error bounds for eigenfunctions}, 
  arXiv preprint arXiv:1904.07903v1, (2019).
  
  \bibitem{MehMie2011}
  Volker Mehrmann and Agnieszka Miedlar, \emph{Adaptive computation of smallest
    eigenvalues of self-adjoint elliptic partial differential equations}, Numer.
    Linear Algebra Appl. \textbf{18} (2011), no.~3, 387--409. \MR{2760060}
  
  \bibitem{Meyer:2000}
  Carl Meyer, \emph{Matrix analysis and applied linear algebra}, Society for
    Industrial and Applied Mathematics (SIAM), Philadelphia, PA, 2000.
    \MR{1777382}
  
  \bibitem{moore2009introduction}
  Ramon~E Moore, R~Baker Kearfott, and Michael~J Cloud, \emph{Introduction to
    interval analysis}, vol. 110, {SIAM}, 2009.
  
  \bibitem{Nakatsukasa2020}
  Yuji Nakatsukasa, \emph{{Sharp error bounds for Ritz vectors and approximate
    singular vectors}}, Math. Comp. \textbf{89} (2020), 1843--1866.
  
  \bibitem{PraSyn:1947}
    William Prager and John Lighton Synge, \emph{Approximations in elasticity based on the concept of function space}, Quart. Appl. Math. \textbf{5} (1947), 241--269.

  \bibitem{Rep-2008-book}
  Sergey Repin,
  \emph{A posteriori estimates for partial differential equations}, 
  Walter de Gruyter GmbH {\&} Co. KG, Berlin, 2008.

  \bibitem{SebVej2014}
  Ivana {\v{S}}ebestov{\'a} and Tom{\'a}{\v{s}} Vejchodsk{\'y}, \emph{Two-sided
    bounds for eigenvalues of differential operators with applications to
    {F}riedrichs, {P}oincar\'e, trace, and similar constants}, SIAM J. Numer.
    Anal. \textbf{52} (2014), no.~1, 308--329. \MR{3163245}
  
  \bibitem{toyonaga2002verified}
  Kenji Toyonaga, Mitsuhiro~T Nakao, and Yoshitaka Watanabe, \emph{Verified
    numerical computations for multiple and nearly multiple eigenvalues of
    elliptic operators}, Journal of Computational and Applied Mathematics
    \textbf{147} (2002), no.~1, 175--190.
  
  \bibitem{Vejchodsky2018b}
  Tom{\'a}{\v{s}} Vejchodsk{\'y}, \emph{Three methods for two-sided bounds of
    eigenvalues--a comparison}, Numer. Methods Partial Differential Equations
    \textbf{34} (2018), no.~4, 1188--1208.
  
  \bibitem{Vejchodsky2018}
  Tom\'a\v{s} Vejchodsk\'y, \emph{Flux reconstructions in the
    {L}ehmann-{G}oerisch method for lower bounds on eigenvalues}, J. Comput.
    Appl. Math. \textbf{340} (2018), 676--690. \MR{3807831}
  
  \bibitem{Xie2LIU-2018}
  Manting Xie, Hehu Xie, and Xuefeng Liu, \emph{{Explicit lower bounds for Stokes
    eigenvalue problems by using nonconforming finite elements}}, Jpn. J. Ind.
    Appl. Math. \textbf{35} (2018), no.~1, 335--354.
  
  \bibitem{Yang2010}
  Yidu Yang, Zhimin Zhang, and Fubiao Lin, \emph{Eigenvalue approximation from
    below using non-conforming finite elements}, Science in China Series A:
    Mathematics \textbf{53} (2010), no.~1, 137--150.
  
  \bibitem{you-xie-liu-2019}
  Chun'guang You, Hehu Xie, and Xuefeng Liu, \emph{Guaranteed eigenvalue bounds
    for the Steklov eigenvalue problem}, SIAM J. Numer. Anal. \textbf{57} (2019),
    no.~3, 1395--1410.
  
\bibitem{Liu-2022}
    Xuefeng Liu, Mitsuhiro T. Nakao, and Shin’ichi Oishi,
    Computer-assisted proof for the stationary solution existence of the Navier–Stokes equation over 3D domains,
    Communications in Nonlinear Science and Numerical Simulation,
    Volume 108,
    2022,
    106223.

 
\end{thebibliography}

\end{document}